\documentclass[twocolumn]{autart}    % Enable this line and disable the
\usepackage{wrapfig}
\usepackage{graphicx}          % Include this line if your
\usepackage{epstopdf}
\usepackage{mathrsfs}
\usepackage{amssymb}
\usepackage{amsmath}
\usepackage{amsfonts}
\usepackage{bbm}
\usepackage{mathrsfs}
\usepackage{subfigure}
\usepackage{float}
\usepackage{color}
\begin{document}

\begin{frontmatter}
\title{Linear-Quadratic Mean-Field Game for Stochastic Systems with Partial Observation} 

\thanks[footnoteinfo]{Corresponding author: Na Li.}
\author[CCE]{Min Li}\ead{lim@sdu.edu.cn},    % Add the
\author[CCD]{Na Li\thanksref{footnoteinfo}}\ead{lina@sdufe.edu.cn},
\author[CCC]{Zhen Wu}\ead{wuzhen@sdu.edu.cn}            % e-mail address
%\author[Baiae]{Publius Maro Vergilius}\ead{vergilius@culture.ir}  % (ead) as shown

\address[CCE]{Institute of Geotechnical and Underground Engineering, Shandong
University, Jinan 250061, China}
\address[CCD]{School of Statistics and Mathematics, Shandong University of Finance and Economics, Jinan 250014,
China}
\address[CCC]{School of Mathematics, Shandong
University, Jinan 250100, China}

 % Please supply
%\address[Rome]{}             % full addresses
%\address[Baiae]{The White House, Baiae}        % here.

\begin{keyword}
Linear-quadratic mean-field game; Partial information; Forward-backward stochastic differential equation;  Optimal filtering; $\varepsilon$-Nash equilibrium.

                       % Five to ten keywords,        % chosen from the IFAC
\end{keyword}                             % keyword list or with the
                                          % help of the Automatica
                                          % keyword wizard

\begin{abstract}   % Abstract of not more than 200 words.
This paper is concerned with a class of linear-quadratic stochastic large-population problems with partial information, where the individual agent only has access to a noisy observation process related to the state. The dynamics of each agent follows a linear stochastic differential equation driven by individual noise, and all agents are coupled together via the control average term. Using the mean-field game approach and the backward separation principle with a state decomposition technique, the decentralized optimal control can be obtained in the open-loop form through a forward-backward stochastic differential equation with the conditional expectation. The optimal filtering equation is also provided. By the decoupling method, the decentralized optimal control can also be further presented as the feedback of state filtering via the Riccati equation.  The explicit solution of the control average limit is given, and the consistency condition system is discussed. Moreover, the related $\varepsilon$-Nash equilibrium property is verified. To illustrate the good performance of theoretical results, an example in finance is studied.
\end{abstract}

\end{frontmatter}

\section{Introduction}
Recently, researchers have made great contributions to developing large-population systems, which can normally be viewed as complicated systems involving considerable interrelated agents. The distinctive feature of large-population problems is that the influence of the individual agent is negligible and insignificant from a micro perspective, while the collective behavior of all agents is significant and cannot be ignored from a macro perspective. Mathematically, there exists a weak coupling relation among the considerable agents in their dynamics or cost functionals through the state average or more general empirical measure. With its own special structure, the large-population problem has drawn increasing attention in many fields, such as engineering, economics, finance, social science, operations research,  management, and so on.

However, when the number of agents is sufficiently large, this kind of problem is intractable to deal with directly due to the curse of dimensionality and the weak coupling structure. For these reasons, it is not unrealistic for an individual agent to obtain a centralized strategy, which also needs the information of other peers'. An alternative solution is to design the so-called decentralized control based on the local information, where only information of the given individual agent and some off-line quantities are enough.  It should be emphasized that the key role in achieving the above goals is the mean-field game (MFG) approach. For a comprehensive understanding of MFG theory, interested readers can refer to \cite{Lions,Cardaliaguet2010,Bensoussan2013,Carmona2018,Carmona20182}.

The theory of MFG was originally formulated by Larsy, Lions \cite{Lasry2007} and Huang, Malham\'{e}, Caines \cite{Huang2006} independently. Compared with the nonlinear counterpart, the dynamic optimization problems for linear stochastic large-population systems enjoy the favor. In particular, the linear-quadratic (LQ) case has been studied extensively for its elegant properties and broad applications. Motivated by some phenomena in engineering, \cite{Huang2007} put forward the Nash certainty equivalence (NCE) principle to study a class of LQ Gaussian (LQG) stochastic large-population problems. The authors in \cite{Li2008} and \cite{Bardi2012} considered the stochastic large-population problem in the infinite horizon and discussed ergodic MFGs. By adopting the adjoint equation approach, \cite{Bensounssan2016} provided a comprehensive study of LQG MFG, where the optimal mean field term satisfied a forward-backward ordinary differential equation. The more general cases, where the agent was subjected to different modelings of state, have also been studied,  see \cite{Wang2012} the Markov jump parameters case, see \cite{Huang2018} the time-delay case, \cite{Li2020} the Poisson jumps case and \cite{Li2022} the non-monotone case. For the mixed MFGs, the readers can see \cite{Huang2010,Nguyen2012,Caines2013,Nie2018}. The cooperative game among the decision makers, known as social optimization, has also been further investigated, see \cite{Malhame2012,Wang2020,Huang2021}.

Generally, the full information is not available in practice, which means that the agents only have access to partial information. Because of its potential theory and application values, the partial information problem has attracted researchers' intensive attention. When dealing with such kinds of problems, the filtering technique is an indispensable and valid tool. The well-known results of linear filtering theory were founded by \cite{Kalman1961}. For nonlinear cases, the celebrated monograph  \cite{Liptser1977} provided a perfect response. The recent result can also be seen \cite{Xiong2008}. To solve the stochastic optimal control problem under partial information, \cite{Wonham1968} proposed the classical separation principle. It suggested that this problem can be treated in a separate manner by first studying the filtering problem and then solving the equivalent optimal control problem.

Inspired by the above observations and the idea in \cite{Bensoussan1992}, a new-style separation principle named backward separation principle was originated by \cite{Wang2008} to solve the partially observable stochastic LQ problem. Here, the component ``backward'' indicates that for research on the combined problem of control and filtering, we first consider the optimal control problem and then study its filtering issue, where the order of processing the problem is contrary to the aforementioned results in \cite{Wonham1968}. There exist some associated works following the backward separation methods. \cite{Huang2009} studied the partial information control problem of backward stochastic systems and derived a new stochastic maximum principle. On the strength of Girsanov's theorem, \cite{Wang2009} and \cite{Wu2010} were concerned with the optimal control problem for partially observed forward-backward stochastic systems, and the corresponding optimal controls were constructed in the sense of a weak solution. By introducing the backward separation principle with a state decomposition technique, \cite{Wang2015} discussed a partial information LQ optimal control problem derived by forward-backward stochastic differential equations (FBSDEs). \cite{Na2020} studied an LQ optimal control problem for time-delay stochastic systems with recursive utility under full and partial information. For the systematic introduction about the partial information optimal control problems of FBSDEs, interested readers can see the monograph \cite{Wu2018}.

For the literature review, some works related to partial information MFG problems should be described. To our best knowledge, the earliest research on this issue can be traced back to \cite{Huang2006meeting}, which studied a partial observation distributed decision-making problem in a system of uniform coupled agents, and each agent only has local noisy measurements of its own state. \cite{Wang2016} explored a class of partial information LQ large-population problems, where the information available to the individual agent was the subset of the whole filtration. Distinguished from the listed work, \cite{Wu2016} originally considered the dynamic optimization problem of backward stochastic large-population systems under partial information. Recently, there also exist some works on partially observed MFGs, see e.g. \cite{Caines2016,Caines2019,Firoozi2020} for major-minor cases, \cite{Bensoussan2019} and \cite{LiNieWu2022} for LQ cases.

Concerned with the dynamic optimization problem of linear stochastic large-population systems, most existing literature focused on the state average type weak coupling structure. In fact, the control average large-population problem is also well documented in MFG theory, see \cite{Huang2018,Li2020,Li2022} and \cite{Huang2013} for more details. Different from the traditional large-population problem setting, the introduction of a control average will bring some new features. On the one hand, the control average can better reflect the instantaneous and immediate effects in large-population systems, which is motivated by some practical applications. For example, we consider the horizontal competition problem in economics. When one company's policy changes, the other companies will make corresponding adjustments. Thus, this kind of problem can be deduced to the control average large-population problem. On the other hand, by comparison, the introduction of a control average will make the definition of the NCE or consistency condition (CC) system different, and the corresponding processing techniques need to be adjusted.

Motivated by the above phenomena, this paper studies a class of control average stochastic large-population problems with partial observation. The main \emph{contributions} of this paper can be summarized as follows.
\begin{itemize}
  \item A class of partial information LQ stochastic large-population problems are considered, where the individual agent only has access to a noisy observation process related to the state. Distinguished from most existing literature, all agents are coupled together through the \emph{control average} structure in our paper.
  \item To break the circular dependence between control and observation in the current setup, the MFG approach and the backward separation principle with a state decomposition technique are adopted to design the decentralized optimal control. Drawing support from an FBSDE involving with the conditional expectation, the decentralized optimal control can be expressed in the open-loop form. The optimal filtering equation is also provided. By the decoupling method, the decentralized optimal control can also be presented as the feedback of state filtering. The CC system is discussed, where the control average limit can be solved explicitly. Moreover, the related $\varepsilon$-Nash equilibrium property is verified.
  \item As an application, a cash management problem with partial observation is well investigated, where the explicit decentralized optimal capital injection or withdrawal is obtained. Some simulation results are also attached.
\end{itemize}

 In the large-population problem of this paper, we focus on the coupling of controls in both dynamics and costs because the coupling via controls of the other agents has \emph{instantaneous} and \emph{immediate} impact on the system. When the coupling of states is incorporated into both dynamics and costs, the results in the paper can be similarly extended. To avoid heavy notation, only the control coupling is introduced and considered here.

The rest of this paper is organized as follows. Section 2 includes three subsections. In subsection 2.1, we present some preliminaries and formulate the LQ large-population problem with partial observation. We study the limiting problem with partial observation and the corresponding MFG in subsection 2.2. Using the backward separation with a state decomposition technique, we derive the decentralized optimal control in the open-loop form. The corresponding optimal filtering equation is established. We also obtain the decentralized optimal control in the feedback form of state filtering. In subsection 2.3, we get the explicit control average limit and discuss the corresponding CC system. Section 3 verifies the related $\varepsilon$-Nash equilibrium property. A practical example is presented in Section 4 to illustrate the effectiveness of theoretical results, and some numerical simulations are attached.  Section 5 concludes the paper.

\section{LQ MFG with Partial Observation}
\subsection{Preliminaries and Problem Formulation}
For a fixed $T>0$, we consider a large-population problem involving $N\in\mathbb{N}$ agents on finite time horizon $[0,T]$, where $\mathbb{N}$ denotes the set of all nonnegative integers. Let $(\Omega,\mathcal{F},\mathbb{F},\mathbb{P})$ be a complete filtered probability space with usual conditions, on which are defined two independent $nN$-dimensional standard Brownian motions $\{W_i(t),1\leq i \leq N \}_{0\leq t \leq T}$ and $\{\overline{W}_i(t),1\leq i \leq N \}_{0\leq t \leq T}$. The main reasons that we need complete filtered space here are listed as follows: (i) The completeness of the probability space guarantees the existence of a continuous modified version for the stochastic process, which is the basis of the well-posedness for some equations. (ii) The stochastic integrals need to be defined on a complete filtered probability space. For $1\leq i \leq N$, $W_i$ denotes the individual noise only for the $i$th agent and let $\mathcal{F}_t^i$ be the filtration generated by $\{W_i(s),0\leq s \leq t\}$. Suppose $\mathcal{F}_t$ be the natural filtration generated by $\{W_i(s),\overline{W}_i(s),0\leq s \leq t,1\leq i \leq N \}$ augmented by all $\mathbb{P}$-null sets $\mathcal{N}$ and $\mathbb{F}=\{\mathcal{F}_t\}_{0\leq t \leq T}$.

Throughout the paper, let $\mathbb{R}^n$ be an $n$-dimensional Euclidean space with usual norm $|\cdot|$ and inner product $\langle\cdot,\cdot\rangle$. $\mathbb{R}^{l\times n}$ denotes $(l\times n)$-dimensional matrix, $I$ represents the identity matrix with appropriate dimension. Let $\mathcal{S}^l$ be the set of all $(l\times l)$-dimensional symmetric matrices. If $\mathcal{M}\in \mathcal{S}^l$ is positive (semi-positive) definite, we write $\mathcal{M}>(\geq)~0$. For any vector or matrix $A$, we use $A^{\top}$ to represent its transpose. For any Euclidean space $\mathcal{H}$ and filtration $\mathcal{G}$, if $g(t,\omega):[0,T]\times\Omega\rightarrow\mathcal{H}$ is a $\mathcal{G}_t$-progressively measurable stochastic process such that $\mathbb{E}\int_0^T|g(t)|^2dt<\infty$, we write $g(t)\in L^2_{\mathcal{G}}(0,T;\mathcal{H})$; if $g(t,\omega):[0,T]\times\Omega\rightarrow\mathcal{H}$ is a continuous $\mathcal{G}_t$-progressively measurable stochastic process such that $\mathbb{E}[\sup_{0\leq t \leq T}|g(t)|^2]<\infty$, we write $g(t)\in \mathcal{L}^2_{\mathcal{G}}(0,T;\mathcal{H})$; if $g(t):[0,T]\rightarrow\mathcal{H}$ is a uniformly bounded deterministic function, we write $g(t)\in L^\infty(0,T;\mathcal{H})$; if $g(t):[0,T]\rightarrow\mathcal{H}$ is a deterministic function such that $ \int_0^T |g(t)|^2dt<\infty$, we write $g(t)\in L^2(0,T;\mathcal{H})$; if $\xi:\Omega\rightarrow\mathcal{H}$ is a $\mathcal{G}_T$-measurable random variable such that $\mathbb{E}[\xi^2]<\infty$, we write $\xi\in L^2_{\mathcal{G}_T}(\Omega,\mathcal{H})$.

Now, we assume that the dynamics of the $i$th agent satisfies the weakly coupled SDE as follows:
\begin{equation}
\left\{
\begin{aligned}
dX_{i}(t)=&~[A(t)X_{i}(t)+B(t)u_{i}(t)+\widetilde{B}(t)u^{(N,-i)}(t)]dt \\
&+\sigma (t)dW_{i}(t), \\
X_{i}(0)=&~x,\label{state}
\end{aligned}
\right.
\end{equation}
where $x\in\mathbb{R}^n$ is the initial value, $X_i(\cdot)$, $u_i(\cdot)$ denote state process and control process of the agent $i$, respectively. Here, we define the control average as
\[u^{(N,-i)}(\cdot):=\frac{1}{N-1}\underset{j=1,j\neq i}{\overset{N}\sum} u_j(\cdot).   \]
For $1\leq i \leq N$, the agent $i$ cannot observe the real state $X_i(\cdot)$ due to some external factors and its own technical limitations. We suppose that each agent can only have access to a noisy observation process related to the state via
\begin{equation}
\left\{
\begin{aligned}\label{observation}
dV_{i}(t)=&~[F(t)X_{i}(t)+G(t)]dt+H(t)d\overline{W}_i(t), \\
V_{i}(0)=&~0.
\end{aligned}%
\right.
\end{equation}
\\
Here, $A(\cdot),B(\cdot),\widetilde{B}(\cdot),\sigma(\cdot)$ $F(\cdot),G(\cdot),H(\cdot)$ are the deterministic functions with appropriate dimensions, and the specific assumptions will be given later.

For simplicity, let $\mathcal{F}^{V_i}_t$ be the augmentation of $\sigma\{V_i(s), $\\$0\leq s \leq t\}$, which stands for the observable information of the $i$th agent up to time $t$. The admissible control set is defined by
$\mathcal{U}_{i}=\{u_{i}(\cdot )|u_{i}(\cdot )\in L_{\mathcal{F}^{V_i}}^{2}(0,T;%
\mathbb{R}^{k})\}$.

Moreover, we use $\boldsymbol{u}(\cdot)=(u_1(\cdot),\ldots,u_N(\cdot))$ to denote the strategy set of all agents and $\boldsymbol{u}_{-i}(\cdot)=(u_1(\cdot),\ldots,u_{i-1}(\cdot),u_{i+1}(\cdot),\ldots,u_N(\cdot))$ to denote the strategy set except the agent $i$, $1\leq i\leq N$. The cost functional is supposed to be
\begin{equation}
\begin{aligned}
&\mathcal{J}_{i}(u_{i}(\cdot ),\boldsymbol{u}_{-i}(\cdot ))=\frac{1}{2}\mathbb{E}\Big\{%
\int_{0}^{T}[\langle Q(t)X_{i}(t),X_{i}(t)\rangle\\
&+\langle
R(t)(u_{i}(t)-K u^{(N,-i)}(t)),u_{i}(t)-K u^{(N,-i)}(t)\rangle ]dt\\
& +\langle MX_{i}(T),X_{i}(T)\rangle \Big\},\label{cost}
\end{aligned}
\end{equation}
where the second term penalizes the departure from the control average.

For the coefficients of \eqref{state}-\eqref{cost}, we give the following assumptions.\\
(A1) The coefficients of the state system \eqref{state} and the observation system \eqref{observation} satisfy
\begin{equation*}
\left\{
\begin{aligned}
&A(\cdot), F(\cdot)\in L^{\infty}(0,T;\mathbb{R}^{n\times n}), \\
&B(\cdot),\widetilde{B}(\cdot)\in L^{\infty}(0,T;\mathbb{R}^{n\times k}),\\
&\sigma(\cdot), G(\cdot),H(\cdot)\in L^{\infty}(0,T;\mathbb{R}^{n\times n}).
\end{aligned}
\right.
\end{equation*}\\
(A2) The coefficients of the cost functional \eqref{cost} satisfy
\begin{equation*}
\left\{
\begin{aligned}
&Q(\cdot)\in L^{\infty}(0,T;\mathcal{S}^{n}), R(\cdot)\in L^{\infty}(0,T;\mathcal{S}^{k}),\\
&K\in\mathcal{S}^k,M\in\mathcal{S}^n, Q(\cdot)\geq 0 , R(\cdot)>0, M\geq 0,\\
&(I-K)~\text{is an invertible matrix.}
\end{aligned}
\right.
\end{equation*}
Under the above assumptions, for any $u_i(\cdot)\in \mathcal{U}_i$, the state system \eqref{state} and the observation system \eqref{observation} admit unique solutions (see \cite[Theorem 6.3 of Chapter 1]{Yongzhou1999}), then the cost functional \eqref{cost} is well-defined. Now, we formulate the dynamic optimization problem for large-population systems with partial observation (LPO) as follows.

\textbf{Problem (LPO)} To find a strategy set $\boldsymbol{u}^{\ast}(\cdot)=(u^{\ast}_1(\cdot),\ldots,u^{\ast}_N(\cdot))$, where $u^{\ast}_i(\cdot)\in \mathcal{U}_i$, such that
\begin{equation*}
\mathcal{J}_{i}(u^{\ast}_{i}(\cdot ),\boldsymbol{u}^{\ast}_{-i}(\cdot ))=\underset{u_i(\cdot)\in\mathcal{U}_i}{\inf}
\mathcal{J}_{i}(u_{i}(\cdot ),\boldsymbol{u}^{\ast}_{-i}(\cdot )),~~1\leq i \leq N,
\end{equation*}
subjects to \eqref{state}-\eqref{cost}.

\subsection{The Limiting Problem with Partial Observation}
It is mentioned that although there exists a finite number of agents in this setting, the complicated coupling mechanism and the curse of dimensionality make it unsolvable. As a matter of fact, it is infeasible for each agent to obtain the centralized strategies in the non-cooperative game framework. Alternatively, by the MFG approach, we are devoted to determining the approximate equilibrium. Roughly speaking, we derive the decentralized control when agent number $N$ tends to infinity and verify that they satisfy $\varepsilon$-Nash equilibrium property.

We note that the key step in MFG is to introduce a frozen control average limit, which results in the so-called limiting problem. When $N\rightarrow\infty$, suppose that $\frac{1}{N-1}\sum_{j=1,j\neq i}^N u_j(\cdot)$ is approximated by $m(\cdot)$. It follows that $m(\cdot)$ is a deterministic function belonging to $\in L^2(0,T;\mathbb{R}^k)$, which will be further verified in later analysis.  We introduce the following auxiliary system
\begin{equation}
\left\{
\begin{aligned}
dX_{i}(t)=&~[A(t)X_{i}(t)+B(t)u_{i}(t)+\widetilde{B}(t)m(t)]dt \\
&+\sigma (t)dW_{i}(t), \\
X_{i}(0)=&~x.\label{lstate}
\end{aligned}
\right.
\end{equation}
The corresponding observation process $V_{i}(\cdot)$ still satisfies \eqref{observation} in form but with $X_{i}(\cdot)$ given by \eqref{lstate}. Moreover, we denote $\mathcal{F}^{V_i}$ with respect to $V_i(\cdot)$ related to $X_i(\cdot)$ of (4) in the following part.

Now, we will specify the conditions that admissible control satisfies. It is natural to determine the control $u_i(\cdot)$ by the observation $V_i(\cdot)$. However, $V_i(\cdot)$ depends on $u_i(\cdot)$ via the state $X_i(\cdot)$, which means that $u_i(\cdot)$ also has the impact on $V_i(\cdot)$.  It implies that there exists \emph{circular dependence} between control and observation in the limiting problem, which brings some difficulties in determining the optimal control. To overcome this issue, we first introduce the following processes
\begin{equation}
\left\{
\begin{aligned}\label{tildeX}
d\tilde{X}_{i}(t)=&~A(t)\tilde{X}_{i}(t)dt+\sigma (t)dW_{i}(t), \\
\tilde{X}_{i}(0)=&~x,%
\end{aligned}
\right.
\end{equation}
and
\begin{equation}
\left\{
\begin{aligned}\label{tildeV}
d\tilde{V}_{i}(t)=&~F(t)\tilde{X}_{i}(t)dt+H(t)d\overline{W}_{i}(t), \\
\tilde{V}_{i}(0)=&~0.%
\end{aligned}
\right.
\end{equation}

For any $u_i(\cdot)\in L_{\mathcal{F}^{V_i}}^{2}(0,T;\mathbb{R}^{k})$ and given $m(\cdot)$, let $\breve{X}_{i}(\cdot)$ and $\breve{V}_{i}(\cdot)$ solve
\begin{equation}
\left\{
\begin{aligned}\label{breveX}
d\breve{X}_{i}(t)=&~[A(t)\breve{X}_{i}(t)+B(t)u_{i}(t)+\widetilde{B}(t)m(t)]dt, \\
\breve{X}_{i}(0)=&~0,%
\end{aligned}
\right.
\end{equation}
and
\begin{equation}
\left\{
\begin{aligned}\label{breveV}
d\breve{V}_{i}(t)=&~[F(t)\breve{X}_{i}(t)+G(t)]dt, \\
\breve{V}_{i}(0)=&~0,
\end{aligned}
\right.
\end{equation}
respectively. If we define
\begin{equation*}
X_{i}(t):=\tilde{X}_{i}(t)+\breve{X}_{i}(t),~V_{i}(t):=\tilde{V}_{i}(t)+\breve{V%
}_{i}(t),
\end{equation*}
 by a simple addition of \eqref{tildeX}-\eqref{breveV}, it can verify that $X_i(\cdot)$ and $V_i(\cdot)$ satisfy \eqref{lstate}, and the corresponding observation process, respectively.

Let $\mathcal{F}_t^{\tilde{V}_i}=\sigma\{\tilde{V}_i(s), 0\leq s \leq t  \} $. The decentralized control set is given by $\mathcal{U}_i^{d}=\{u_{i}(\cdot )|u_{i}(\cdot )\in\widetilde{\mathcal{U}_i}$ is $\mathcal{F}^{V_i}$-adapted$\}$,
where $\widetilde{\mathcal{U}_i}=\{u_i(\cdot)|u_i(\cdot)$ is an $\mathcal{F}^{\tilde{V}_i}$-adapted process in $\mathbb{R}^k$ satisfying $\mathbb{E}[\underset{0\leq t \leq T} {\sup}u_i(t)^2]<\infty  \}$.

 \textbf{Remark 2.1.} It is worth stressing that the decentralized control here should be adapted to both $\mathcal{F}^{V_i}$ and $\mathcal{F}^{\tilde{V}_i}$. The main reasons are listed as follows. On the one hand, it is natural to define that the decentralized control $u_i(\cdot)\in L^2_{\mathcal{F}^{V_i}}(0,T;\mathbb{R}^k)$ in the partially observed large-population problem. On the other hand, to break the \emph{circular dependence} between control and observation in the limiting problem, we split the state/observation process of the $i$th agent into two parts: ($\tilde{X}_i(\cdot),\tilde{V}_i(\cdot)$) and ($\breve{X}_i(\cdot),\breve{V}_i(\cdot)$), where the first part is independent of control. Moreover, we will prove that in such a decentralized control set, filtrations $\mathcal{F}^{V_i}$ and $\mathcal{F}^{\tilde{V}_i}$ are consistent (see Lemma 2.1).

Then the limiting cost functional becomes
\begin{equation}
\begin{aligned}
J_{i}(u_{i}(\cdot ))
=&\frac{1}{2}\mathbb{E}\Big\{\int_{0}^{T}[\langle
Q(t)X_{i}(t),X_{i}(t)\rangle\\
& +\langle R(t)(u_{i}(t)-Km(t)),u_{i}(t)-Km(t)\rangle
]dt\\
&+\langle MX_{i}(T),X_{i}(T)\rangle \Big\}. \label{lcost}
\end{aligned}
\end{equation}
\\
Now, the original \emph{coupling} problem is decoupled among the agents by fixing $m(\cdot)$ and considering the decentralized control, which turns out to be an LQ stochastic optimal control problem with partial observation. Now, the limiting large-population problem with partial observation (LLPO) can be formulated as follows.

\textbf{Problem (LLPO)} For $i$th agent, find a control $u^{\ast}_i(\cdot)\in \mathcal{U}_i^d$ such that
\begin{equation*}
J_{i}(u^{\ast}_{i}(\cdot ))=\underset{u_{i}(\cdot )\in\mathcal{U}_i^d}{\inf}J_{i}(u_{i}(\cdot )),~~1\leq i \leq N,
\end{equation*}
subjects to \eqref{lstate}, \eqref{observation} and \eqref{lcost}.

When $u^{\ast}_i(\cdot)$ meets above relationship, it is called the decentralized optimal control of Problem (LLPO). The corresponding $X^{\ast}_i(\cdot)$ with respect to $u^{\ast}_i(\cdot)$ is called the decentralized optimal state trajectory, and $(u^{\ast}_i(\cdot),X^{\ast}_i(\cdot))$ is the decentralized optimal pair.

 \textbf{Remark 2.2.} For simplicity, the same notations for the state $X_i(\cdot)$, the control $u_i(\cdot)$ are applied in \eqref{state} and \eqref{lstate} as well as \eqref{cost} and \eqref{lcost}. Since the introduction of the control average limit, the auxiliary state in \eqref{lstate} is quite different from the original state in \eqref{state}. The same arguments also hold true
for the observation process $V_i(\cdot)$.

Before displaying the main results of this section, some useful lemmas are necessary.

\textbf{Lemma 2.1.} For $ i \in \mathbb{N}$ and any $u_i(t)\in \mathcal{U}_i^d$, $\mathcal{F}^{V_i}_t=\mathcal{F}^{\tilde{V}_i}_t$.\\
\emph{Proof.} For any $u_i(t)\in \mathcal{U}_i^d$, since it is $\mathcal{F}^{\tilde{V}_i}_t$-adapted, thus $\breve{X}_i(t)$ is $\mathcal{F}^{\tilde{V}_i}_t$-adapted following from \eqref{breveX}, so is $\breve{V}_i(t)$. From the relationship $V_i(t)=\tilde{V}_i(t)+\breve{V}_i(t)$, we can conclude that $\mathcal{F}^{V_i}_t\subseteq\mathcal{F}^{\tilde{V}_i}_t$. In turn, we can prove that $\mathcal{F}^{\tilde{V}_i}_t\subseteq\mathcal{F}^{V_i}_t$ by the relationship $\tilde{V}_i(t)=V_i(t)-\breve{V}_i(t)$ in the same manner.   \hfill$\square$

Using the above arguments, we introduce a \emph{restrictive} definition of admissibility as  $\mathcal{U}_i^d$ to handle the circular dependence in the limiting problem with partial observation. In the classical literature, this definition and approach have been widely studied, see e.g. \cite{Bensoussan1992}, \cite{Wang2015} and \cite{Bensoussan2019}. In this way, each agent has access to the information augmented by $V_i(\cdot)$ and $\tilde{V}_i(\cdot)$.

We show that the following lemma is applicable to the limiting problem.  Since the proofs are similar to the results in \cite{Wang2015}, we omit the details here.

\textbf{Lemma 2.2.} Let (A1) and (A2) hold. For any given $m(\cdot)\in L^2(0,T;\mathbb{R}^k)$, then
\[\underset{u_{i}(\cdot )\in\mathcal{U}_i^d}{\inf}J_{i}(u_{i}(\cdot ))=\underset{v_{i}(\cdot )\in\widetilde{\mathcal{U}}_i}{\inf}J_{i}(v_{i}(\cdot )), ~~ i \in \mathbb{N}.\]

So far, we have introduced the state decomposition technique to break the circular dependence between control and observation. In the remaining part, we will solve the limiting control problem with partial observation using the backward separation principle. By standard variational calculus and dual method, we first derive the decentralized optimal control for Problem (LLPO) in the open-loop form.  The corresponding optimal filtering equation is established. We also obtain the decentralized optimal control in the feedback form of state filtering. To make this part clearer, we divide the above derivation into three steps as follows.

\textbf{Step 1: The decentralized optimal control in the open-loop form}

In the following theorem, we show that the decentralized optimal control can be expressed in the open-loop form by the FBSDE with conditional expectation. Firstly, we define the open-loop decentralized optimal control as follows. For more details of LQ open-loop optimal control, the readers can refer to \cite[Definition 2.1.3]{SunYong2020}.

\textbf{Definition 2.1.}  For $1\leq i \leq N$, if there exists a unique $u^{\ast}_i(\cdot)\in \mathcal{U}_i^d$ satisfying
\[J_{i}(u^{\ast}_{i}(\cdot ))\leq J_{i}(u_{i}(\cdot )),~~~~\text{for any }  u_{i}(\cdot )\in \mathcal{U}_i^d,\]
 such $u^{\ast}_i(\cdot)$ is called the \emph{open-loop} decentralized optimal control of Problem (LLPO).

\textbf{Theorem 2.1.} Let (A1) and (A2) hold. For any given $m(\cdot)\in L^2(0,T;\mathbb{R}^k)$, suppose $(u^{\ast}_i(\cdot),X^{\ast}_i(\cdot))$ is the decentralized optimal pair of Problem (LLPO), it follows that
\begin{equation}
u^{\ast}_{i}(t)=-R(t)^{-1}B(t)^{\top }\mathbb{E}[Y^{\ast}_{i}(t)|\mathcal{F}%
^{V^{\ast}_i}_t]+Km(t),\label{OC}
\end{equation}
where $(X^{\ast}_{i}(\cdot),Y^{\ast}_{i}(\cdot),Z^{\ast}_{i}(\cdot))\in \mathcal{L}^2_{\mathcal{F}^i}(0,T;\mathbb{R}^n)\times \mathcal{L}^2_{\mathcal{F}^i}(0,T;\mathbb{R}^n)\times L^2_{\mathcal{F}^i}(0,T;\mathbb{R}^n)$ is a solution to the following FBSDE:
\begin{equation}
\left\{
\begin{aligned}\label{FBSDE}
dX^{\ast}_{i}(t)=&~[A(t)X^{\ast}_{i}(t)+B(t)u^{\ast}_{i}(t)+\widetilde{B}(t)m(t)]dt\\
&+\sigma
(t)dW_{i}(t), \\
dY^{\ast}_{i}(t)=&-[A(t)^{\top }Y^{\ast}_{i}(t)+Q(t)X^{\ast}_{i}(t)]dt\\
&+Z^{\ast}%
_{i}(t)dW_{i}(t), \\
X^{\ast}_{i}(0)=&~x,~Y^{\ast}_{i}(T)=MX^{\ast}_{i}(T),%
\end{aligned}
\right.
\end{equation}
and $\mathcal{F}_t^{V^{\ast}_i}=\sigma\{V^{\ast}_i(s),0\leq s \leq t\}$ with $V^{\ast}_i(\cdot)$ satisfies
\begin{equation*}
\left\{
\begin{aligned}
dV^{\ast}_{i}(t)=&~[F(t)X^{\ast}_{i}(t)+G(t)]dt+H(t)d\overline{W}_{i}(t), \\
V^{\ast}_{i}(0)=&~0.
\end{aligned}%
\right.
\end{equation*}
\\
\textbf{Proof.}
Due to the result in Lemma 2.2, we only need to minimize $J_{i}(u_{i}(\cdot ))$ over $\widetilde{\mathcal{U}_i}$.  For any $u_i(\cdot)\in \widetilde{\mathcal{U}}_i$, we have $u_i^\lambda(\cdot):= u_i^\ast(\cdot)+\lambda u_i(\cdot) \in \widetilde{\mathcal{U}}_i$, where $0\leq \lambda \leq  1$. We denote $X_i^\lambda(\cdot)$ is the limiting state corresponding to $u_i^\lambda(\cdot)$, thus it satisfies
\begin{equation*}
\left\{
\begin{aligned}
dX^{\lambda}_{i}(t)
=&~[A(t)X^{\lambda}_{i}(t)+B(t)(u^{\ast}_{i}(t)+\lambda u_i(t))\\
&+\widetilde{B}(t)m(t)]dt+\sigma
(t)dW_{i}(t), \\
X^{\lambda}_{i}(0)=&~x,
\end{aligned}
\right.
\end{equation*}
then
\begin{equation*}
\left\{
\begin{aligned}
&d(X^{\lambda}_{i}(t)-X^{\ast}_{i}(t))\\
=&~[A(t)(X^{\lambda}_{i}(t)-X^{\ast}_{i}(t))+B(t)\lambda u_i(t)]dt,\\
&X^{\lambda}_{i}(0)-X^{\ast}_{i}(0)=~0.
\end{aligned}
\right.
\end{equation*}

By \eqref{lcost} and the inequality \[\langle Qa,a \rangle-\langle Qb,b \rangle\leq \langle Q(a-b),a-b \rangle+2\langle Q(a-b),b \rangle,\] it follows that
$
J_{i}(u_{i}^\lambda(\cdot ))-J_{i}(u_{i}^\ast(\cdot )):
=\Theta_1+\Theta_2,$
where
\begin{equation*}
\begin{aligned}
\Theta_1:=&\frac{1}{2}\mathbb{E}\Big\{\int_{0}^{T}[\langle
Q(t)(X_{i}^\lambda(t)-X_{i}^\ast(t)),X_{i}^\lambda(t)-X_{i}^\ast(t)\rangle\\
& +\langle R(t)(\lambda u_{i}(t)),\lambda u_{i}(t)\rangle
]dt\\
&+\langle M(X_{i}^\lambda(T)-X_{i}^\ast(T)),X_{i}^\lambda(T)-X_{i}^\ast(T)\rangle \Big\}
\end{aligned}
\end{equation*}
and
\begin{equation*}
\begin{aligned}
\Theta_2:=&\mathbb{E}\Big\{\int_{0}^{T}[\langle Q(t)(X_{i}^\lambda(t)-X_{i}^\ast(t)), X_{i}^\ast(t)\rangle\\
&+\langle R(t)(\lambda u_{i}(t)),u_{i}^\ast(t)-Km(t)\rangle]dt\\
&+\langle M(X_{i}^\lambda(T)-X_{i}^\ast(T)),X_{i}^\ast(T)\rangle\Big\}.
\end{aligned}
\end{equation*}

Due to the positive definiteness assumptions of the coefficient matrices, it is obvious that $\Theta_1\geq 0$. Since the optimality of $u_i^\ast(\cdot)$, we need $\Theta_2=0$. To address with $\Theta_2$, we introduce the adjoint process $Y_{i}^\ast(\cdot)$ such that $Y_{i}^\ast(T)=MX_{i}^\ast(T)$ and
\begin{equation*}
dY_{i}^\ast(t)=\Delta(t)dt+Z^{\ast}_{i}(t)dW_{i}(t),
\end{equation*}
where $\Delta(\cdot)$ is an undetermined term. Moreover, $Z^{\ast}_{i}(\cdot)$ is the martingale term which guarantees the adaptness of the solution to the backward equation.

Applying It\^o's formula to $\langle Y_{i}^\ast(\cdot), X^{\lambda}_{i}(\cdot)-X^{\ast}_{i}(\cdot) \rangle$, we have
\begin{equation*}
\begin{aligned}
&d\langle Y_{i}^\ast(t), X^{\lambda}_{i}(t)-X^{\ast}_{i}(t) \rangle\\
=& [ \langle \Delta(t)+ A(t)^{\top }Y^{\ast}_{i}(t),X^{\lambda}_{i}(t)-X^{\ast}_{i}(t)\rangle\\
& + \langle B(t)^{\top}Y^{\ast}_{i}(t), \lambda u_i(t) \rangle    ]dt\\
&+\langle Z^{\ast}_{i}(t)dW_{i}(t), X^{\lambda}_{i}(t)-X^{\ast}_{i}(t)\rangle,
\end{aligned}
\end{equation*}
which implies
\begin{equation*}
\begin{aligned}
&\mathbb{E}\langle M(X_{i}^\lambda(T)-X_{i}^\ast(T)),X_{i}^\ast(T)\rangle\\
=&\mathbb{E}\int_0^T[ \langle \Delta(t)+ A(t)^{\top }Y^{\ast}_{i}(t),X^{\lambda}_{i}(t)-X^{\ast}_{i}(t)\rangle\\
 &+ \langle B(t)^{\top}Y^{\ast}_{i}(t), \lambda u_i(t) \rangle    ]dt.
\end{aligned}
\end{equation*}

Substituting the above equality into $\Theta_2$, we obtain
\begin{equation*}
\begin{aligned}
&\mathbb{E}\int_0^T[ \langle \Delta(t)+ A(t)^{\top }Y^{\ast}_{i}(t)+Q(t)X_{i}^\ast(t),X^{\lambda}_{i}(t)-X^{\ast}_{i}(t)\rangle\\
&+\langle R(t)(u^{\ast}_{i}(t)-Km(t))+B(t)^{\top}Y^{\ast}_{i}(t),\lambda u_i(t)\rangle]dt=0,
\end{aligned}
\end{equation*}
then let
\[\Delta(t)=-[ A(t)^{\top }Y^{\ast}_{i}(t)+Q(t)X_{i}^\ast(t)],\]
 we have
\begin{equation*}
\begin{aligned}
0=&\mathbb{E}\int_0^T\langle R(t)(u^{\ast}_{i}(t)-Km(t))+B(t)^{\top}Y^{\ast}_{i}(t),u_i(t)\rangle dt\\
=&\mathbb{E}\int_0^T\langle R(t)(u^{\ast}_{i}(t)-Km(t))+B(t)^{\top}\mathbb{E}[Y^{\ast}_{i}(t)|\mathcal{F}%
^{\tilde{V}_i}_t],u_i(t)\rangle dt,
\end{aligned}
\end{equation*}
where $Y^{\ast}_{i}(\cdot)$ satisfies the backward component of \eqref{FBSDE}.

We recall  $u^{\ast}_i(\cdot)\in \mathcal{U}_i^d$, then it holds $\mathcal{F}^{\tilde{V}_i}_t=\mathcal{F}^{V_i^{\ast}}_t$ by Lemma 2.1. The above formula yields
\begin{equation*}
R(t)(u^{\ast}_{i}(t)-Km(t))+B(t)^{\top}\mathbb{E}[Y^{\ast}_{i}(t)|\mathcal{F}%
^{V^{\ast}_i}_t]=0,
\end{equation*}
thus the decentralized optimal control can be represented in the open-loop form as \eqref{OC}.  \hfill$\square$

Substituting \eqref{OC} into \eqref{FBSDE}, it can be rewritten as
\begin{equation}
\left\{
\begin{aligned}\label{modFBSDE}
dX^{\ast}_{i}(t)=&~[A(t)X^{\ast}_{i}(t)-B(t)R(t)^{-1}B(t)^{\top }\mathbb{E}[Y^{\ast}_{i}(t)|\mathcal{F}%
^{V^{\ast}_i}_t]\\
&+(B(t)K+\widetilde{B}(t))m(t)]dt+\sigma
(t)dW_{i}(t), \\
dY^{\ast}_{i}(t)=&-[A(t)^{\top }Y^{\ast}_{i}(t)+Q(t)X^{\ast}_{i}(t)]dt+Z^{\ast}%
_{i}(t)dW_{i}(t),\\
X^{\ast}_{i}(0)=&~x,~Y^{\ast}_{i}(T)=MX^{\ast}_{i}(T),
\end{aligned}
\right.
\end{equation}
which is an FBSDE with conditional expectation.

\textbf{Step 2: The optimal filtering}

For $1\leq i\leq N$, if we introduce the notation $\widehat{\theta}_{i}(\cdot ):=\mathbb{E}[\theta_{i}(\cdot )|\mathcal{F}^{%
V^{\ast}_{i}}_\cdot]$. Then the decentralized optimal control can be rewritten as $u^{\ast}_{i}(t)=-R(t)^{-1}B(t)^{\top }\widehat{Y}^{\ast}_{i}(t)+Km(t)$, where $\widehat{Y}^{\ast}_{i}(\cdot)$ satisfies the optimal filtering equation in the following proposition.

\textbf{Proposition 2.1.} Let (A1) and (A2) hold. For any given $m(\cdot)\in L^2(0,T;\mathbb{R}^k)$, the optimal filtering of $(X^{\ast}_{i}(\cdot),Y^{\ast}_{i}(\cdot),Z^{\ast}_{i}(\cdot))$ satisfies
\begin{equation}
\left\{
\begin{aligned}\label{filteringFBSDE}
&d\widehat{X}^{\ast}_{i}(t)\\
=&~[A(t)\widehat{X}^{\ast}_{i}(t)+B(t)u^{\ast}%
_{i}(t)+\widetilde{B}(t)m(t)]dt\\
&+\Pi (t)F(t)^{\top }(H(t)^\top)^{-1}d\widehat{W}_{i}(t),
\\
&d\widehat{Y}^{\ast}_{i}(t)\\
=&-[A(t)^{\top }\widehat{Y}^{\ast}_{i}(t)+Q(t)\widehat{X}^{\ast}_{i}(t)]dt+(\widehat{Y^{\ast}_i(t)X^{\ast}_i(t)^\top}\\
&~~-\widehat{Y}^{\ast}_i(t)\widehat{X}^{\ast}_i(t)^\top) F(t)^\top(H(t)^\top)^{-1}d\widehat{W}_{i}(t), \\
&\widehat{X}^{\ast}_{i}(t)=~x,~\widehat{Y}^{\ast}_{i}(T)=~M\widehat{X}^{\ast}_{i}(T)%
\end{aligned}%
\right.
\end{equation}
with $u^{\ast}_{i}(t)=-R(t)^{-1}B(t)^{\top }\widehat{Y}^{\ast}_{i}(t)+Km(t)$,
where $\widehat{W}_i(\cdot)$ is an $\mathbb{R}^{n}$-valued standard Brownian motion satisfying
\begin{equation}\label{hatW}
\widehat{W}_i(t)= \int_0^t H(s)^{-1}F(s)(X_i^*(s)-\widehat{X}_i^*(s))ds+\overline{W}_i(t)
\end{equation}
and the mean square error $\Pi (\cdot)$ of the estimate $\widehat{X}^{\ast}_{i}(\cdot)$ is the unique solution to
\begin{equation}
\left\{
\begin{aligned}\label{Pi}
&\dot{\Pi}(t)-A(t)\Pi (t)-\Pi (t)A(t)^{\top }\\
&+\Pi
(t)F(t)^\top(H(t)^\top)^{-1}[\Pi (t)F(t)^\top(H(t)^\top)^{-1}]^{\top }\\
& -\sigma (t)\sigma (t)^{\top }=0, \\
&\Pi (0)=0.%
\end{aligned}%
\right.
\end{equation}
\textbf{Proof.} By Lemma 2.1 and $X_{i}^*(t)=\tilde{X}^*_{i}(t)+\breve{X}^*_{i}(t),$ we have
$\widehat{X}^{\ast}_{i}(t)=\mathbb{E}[X^{\ast}_{i}(t)|\mathcal{F}^{V^{\ast}_{i}}_t]=\mathbb{E}[\tilde{X}^{\ast}_{i}(t)|\mathcal{F}^{\tilde{V}_i}_t]
+\breve{X}^{\ast}_{i}(t).$

Applying the multidimensional main filtering equation in Liptser and Shiryaev \cite[Thereom 12.7]{Liptser1977}, we obtain
\begin{equation*}
\begin{aligned}
&\mathbb{E}[X^{\ast}_{i}(t)|\mathcal{F}^{V^{\ast}_{i}}_t]\\
=&~x+\int_0^t\{A(s)^\top\mathbb{E}[X^{\ast}_{i}(s)|\mathcal{F}^{V^{\ast}_{i}}_s]+B(s)u_i^\ast(s)+\widetilde{B}(s)m(s)\}ds\\
&+\int_0^t\Pi_i(s)F(s)^\top(H(s)^\top)^{-1}d\widehat{W}_i(s),
\end{aligned}
\end{equation*}
where for $1\leq i\leq N$, we denote \[\Pi_i (t):=\mathbb{E}[(X^{\ast}_{i}(t)-\widehat{X}^{\ast}_{i}(t))(X^{\ast}_{i}(t)-%
\widehat{X}^{\ast}_{i}(t))^{\top }|\mathcal{F}^{%
V^{\ast}_{i}}_t].\]

 According to the definition of $\Pi_i(\cdot)$, we have
\begin{equation}
\begin{aligned}\label{pirelation}
\Pi_i (t)
=&\mathbb{E}[X^{\ast}_{i}(t)X^{\ast}_{i}(t)^\top-X^{\ast}_{i}(t)\widehat{X}^{\ast}_{i}(t)^\top\\
&-\widehat{X}^{\ast}_{i}(t)X^{\ast}_{i}(t)^\top+\widehat{X}^{\ast}_{i}(t)\widehat{X}^{\ast}_{i}(t)^\top|\mathcal{F}^{%
V^{\ast}_{i}}_t]\\
=&\mathbb{E}[X^{\ast}_{i}(t)X^{\ast}_{i}(t)^\top|\mathcal{F}^{V^{\ast}_{i}}_t]-\widehat{X}^{\ast}_{i}(t)\widehat{X}^{\ast}_{i}(t)^\top\\
=&\widehat{X^{\ast}_{i}(t)X^{\ast}_{i}(t)^\top}-\widehat{X}^{\ast}_{i}(t)\widehat{X}^{\ast}_{i}(t)^\top.
\end{aligned}
\end{equation}

Applying It\^o's formula to $X^{\ast}_{i}(\cdot)X^{\ast}_{i}(\cdot)^\top $ and using \cite[Thereom 12.7]{Liptser1977} again, it follows
\begin{equation}
\begin{aligned}\label{hxx}
&d\widehat{X^{\ast}_{i}(t)X^{\ast}_{i}(t)^\top }\\=&\Big\{A(t)\widehat{X^{\ast}_{i}(t)X^{\ast}_{i}(t)^\top }+\widehat{X^{\ast}_{i}(t)X^{\ast}_{i}(t)^\top } A(t)^\top+B(t)u_i^\ast(t)\widehat{X}^{\ast}_{i}(t)^\top\\
&+\widehat{X}^{\ast}_{i}(t)u_i^\ast(t)^\top B(t)^\top+\widetilde{B}(t)m(t)\widehat{X}^{\ast}_{i}(t)^\top\\
&+ \widehat{X}^{\ast}_{i}(t)m(t)^\top\widetilde{B}(t)^\top+\sigma(t)\sigma(t)^\top\Big\}dt\\
&+\widehat{X^{\ast}_{i}(t)X^{\ast}_{i}(t)^\top X^{\ast}_{i}(t)}(d\widehat{W}_i(t))^\top H(t)^{-1}F(t)\\
&-\widehat{X^{\ast}_{i}(t)X^{\ast}_{i}(t)^\top }F(t)^\top(H(t)^\top)^{-1}d\widehat{W}_i(t)\widehat{X}^{\ast}_{i}(t)^\top.
\end{aligned}
\end{equation}

By the first equation in \eqref{filteringFBSDE} and then applying It\^o's formula to $\widehat{X}^{\ast}_{i}(\cdot)\widehat{X}^{\ast}_{i}(\cdot)^\top $, we have
\begin{equation}
\begin{aligned}\label{hxhx}
&d\widehat{X}^{\ast}_{i}(t)\widehat{X}^{\ast}_{i}(t)^\top \\ =&\Big\{A(t)\widehat{X}^{\ast}_{i}(t)\widehat{X}^{\ast}_{i}(t)^\top+\widehat{X}^{\ast}_{i}(t)\widehat{X}^{\ast}_{i}(t)^\top A(t)^\top \\ &+B(t)u_i^\ast(t)\widehat{X}^{\ast}_{i}(t)^\top+\widehat{X}^{\ast}_{i}(t)u_i^\ast(t)^\top B(t)^\top\\
&+\widetilde{B}(t)m(t)\widehat{X}^{\ast}_{i}(t)^\top+ \widehat{X}^{\ast}_{i}(t)m(t)^\top\widetilde{B}(t)^\top\\
&+\Pi_i(t)F(t)^\top(H(t)^\top)^{-1}[\Pi_i (t)F(t)^\top(H(t)^\top)^{-1}]^{\top }\Big\}dt\\
&+\Pi_i(t)F(t)^\top(H(t)^\top)^{-1}d\widehat{W}_i(t)\widehat{X}^{\ast}_{i}(t)^\top\\
&+\widehat{X}^{\ast}_{i}(t)(d\widehat{W}_i(t))^\top[\Pi_i(t)F(t)^\top(H(t)^\top)^{-1}]^\top.
\end{aligned}
\end{equation}

By taking difference between \eqref{hxx} and \eqref{hxhx} then recalling \eqref{pirelation}, it is easy to verify that $\Pi_i(\cdot)$ independent of $i$ is indeed the solution to \eqref{Pi}, where the proof of the zero diffusion coefficient component is similar to the formula (12.27) in \cite{Liptser1977}. Hence, we omit the index $i$ in $\Pi_i(\cdot)$ and use the notation $\Pi(\cdot)$, which implies the desired results.

As for the backward equation in \eqref{FBSDE}, taking integral from $0$ to $t$, it yields
\begin{equation}
\begin{aligned}
Y^{\ast}_{i}(t)=&Y^{\ast}_{i}(0)-\int_0^t[A(s)^{\top }Y^{\ast}_{i}(s)+Q(s)X^{\ast}_{i}(s)]ds\\
&+\int_0^tZ^{\ast}%
_{i}(s)dW_{i}(s).\label{Pfor}
\end{aligned}
\end{equation}
\\
Using \cite[Theorem 12.7]{Liptser1977} again, we have
\begin{equation*}
\begin{aligned}
&\widehat{Y}^{\ast}_{i}(t)=Y^{\ast}_{i}(0)-\int_0^t[A(s)^{\top }\widehat{Y}^{\ast}_{i}(s)+Q(s)\widehat{X}^{\ast}_{i}(s)]ds\\
&\quad+\int_0^t(\widehat{Y^{\ast}_i(s)X^{\ast}_i(s)^\top}-\widehat{Y}^{\ast}_i(s)\widehat{X}^{\ast}_i(s)^\top)\\
&\quad \times F(s)^\top(H(s)^\top)^{-1}d\widehat{W}_{i}(s).
\end{aligned}
\end{equation*}
\\
Moreover, we have $\widehat{Y}^{\ast}_{i}(T)=M\widehat{X}^{\ast}_{i}(T)$, which implies the objective result.   \hfill$\square$

\textbf{Step 3: The decentralized optimal control in the feedback form}

 Because of the appearance of the conditional expectation term, \eqref{modFBSDE} does not satisfy the standard form like some existing solvable equations in \cite{Ma1999}. Thus, the well-posedness of \eqref{modFBSDE} cannot be guaranteed by the classical FBSDE theory, which hints at the direction of our future work. In the rest of this section, we are going to obtain the feedback representation of decentralized optimal control via the decoupling method.

\textbf{Theorem 2.2.} Let (A1) and (A2) hold. For any given $m(\cdot)\in L^2(0,T;\mathbb{R}^k)$, if $u^{\ast}_i(\cdot)$ is the decentralized optimal control of Problem (LLPO), then it can be expressed as
\begin{equation}
u^{\ast}_i(t)=-R(t)^{-1}B(t)^{\top}(P(t)\widehat{X}^{\ast}_i(t)+\psi(t))+Km(t),\label{statefeedback}
\end{equation}
where
\begin{equation}
\left\{
\begin{aligned}\label{Theta}
&\dot{P}(t)+P(t)A(t)+A(t)^{\top}P(t)\\
&\qquad-P(t)B(t)R(t)^{-1}B(t)^{\top}P(t) +Q(t)=0,\\
&P(T)=M,
\end{aligned}
\right.
\end{equation}
and
\begin{equation}
\left\{
\begin{aligned}\label{psi}
&\dot{\psi}(t)+(A(t)^{\top}-P(t)B(t)R(t)^{-1}B(t)^{\top})\psi(t)\\
&\qquad +P(t)(B(t)K+\widetilde{B}(t))m(t)=0,\\
&\psi(T)=0.
\end{aligned}
\right.
\end{equation}
Moreover, the decentralized optimal state filtering is the solution to
\begin{equation}
\left\{
\begin{aligned}\label{closedstate}
& d\widehat{X}^{\ast}_{i}(t)\\
=&~[(A(t)-B(t)R(t)^{-1}B(t)^{\top}P(t))\widehat{X}^{\ast}_{i}(t)\\
&-B(t)R(t)^{-1}B(t)^{\top}\psi(t)+(B(t)K+\widetilde{B}(t))m(t)]dt\\
&+\Pi (t)F(t)^{\top }(H(t)^\top)^{-1}d\widehat{W}_{i}(t),\\
&\widehat{X}^{\ast}_{i}(0)=~x
\end{aligned}
\right.
\end{equation}
with $\widehat{W}_{i}(\cdot)$, $\Pi(\cdot)$ given by \eqref{hatW} and \eqref{Pi}, respectively.

\textbf{Proof.} According to the terminal condition of \eqref{filteringFBSDE}, we conjecture that
\begin{equation}\label{hatP}
\widehat{Y}^{\ast}_i(t)=P(t)\widehat{X}^{\ast}_i(t)+\psi(t)
\end{equation}
with $P(T)=M$ and $\psi(T)=0$. Substituting \eqref{hatP} into \eqref{OC}, the feedback representation of the decentralized optimal control is given by \eqref{statefeedback}.

Applying It\^o's formula to $\widehat{Y}^{\ast}_i(\cdot)$, it follows that
\begin{equation}
\begin{aligned}\label{YITO}
d\widehat{Y}^{\ast}_i(t)=&[(\dot{P}(t)+P(t)A(t))\widehat{X}^{\ast}_i(t)+P(t)B(t)u^{\ast}_i(t)\\
&+P(t)(B(t)K+\widetilde{B}(t))m(t)+\dot{\psi}(t)]dt\\
&+P(t)\Pi (t)F(t)^{\top }(H(t)^\top)^{-1}d\widehat{W}_{i}(t).
\end{aligned}
\end{equation}

Comparing \eqref{YITO} with the drift term of the second equation in \eqref{filteringFBSDE}, then by noting \eqref{statefeedback} and \eqref{hatP}, one gets
\begin{equation*}
\begin{aligned}
&(\dot{P}(t)+P(t)A(t)+A(t)^{\top}P(t)\\
&-P(t)B(t)R(t)^{-1}B(t)^{\top}P(t)+Q(t))\widehat{X}^{\ast}_i(t)\\
&+\dot{\psi}(t)+(A(t)^{\top}-P(t)B(t)R(t)^{-1}B(t)^{\top})\psi(t)\\
&+P(t)(B(t)K+\widetilde{B}(t))m(t)=0,
\end{aligned}
\end{equation*}
which naturally implies the equations of \eqref{Theta} and \eqref{psi}.

Moreover, substituting \eqref{statefeedback} into the first equation of \eqref{filteringFBSDE}, we can prove that $\widehat{X}^{\ast}_{i}(\cdot)$ is the solution to \eqref{closedstate}, which completes the proof.
   \hfill$\square$

\textbf{Remark 2.3.} By \cite[Theorem 12.7]{Liptser1977} and the equation for $V_i^*(\cdot)$, we have $\widehat{W}_i(\cdot)$ satisfying \eqref{hatW}. After some calculations, we obtain that the equation for $X_i^*(\cdot)-\widehat{X}_i^*(\cdot)$ is driven by $W_i(\cdot)$ and $\overline{W}_i(\cdot)$. Since $W_i(\cdot)$ and $\overline{W}_i(\cdot)$ are both independent identical distribution (i.i.d.), and $P(\cdot)$, $\psi(\cdot)$, $m(\cdot)$ are deterministic, we have that $\widehat{X}_i^*(\cdot)$ is i.i.d. from \eqref{hatW} and \eqref{closedstate}, thus $\widehat{Y}_i^*(\cdot)$ is i.i.d. from \eqref{hatP}.  Moreover,  if the diffusion term is dependent of state and control processes, the backward separation principle and the partial observable results in our paper do not hold, thus $\widehat{Y}_i^*(\cdot)$ is not i.i.d..

\subsection{CC System }
Up to now, all the analysis are established on the fact that $m(\cdot)$ is given. Next, the CC system is introduced to characterize a structure of the control average limit.

\textbf{Proposition 2.2.} Let (A1) and (A2) hold. Then the control average limit $m(\cdot)\in L^2(0,T;\mathbb{R}^k)$ satisfying
\begin{equation}\label{m}
m(t)=-(I-K)^{-1}R(t)^{-1}B(t)^{\top}(P(t)X(t)+\psi(t)),
\end{equation}
where $(X(\cdot),\psi(\cdot))$ is determined by
\begin{equation}
\left\{
\begin{aligned}\label{CC}
&dX(t)\\=&~\{[A(t)-B(t)R(t)^{-1}B(t)^{\top}P(t)\\
&-(B(t)K+\widetilde{B}(t))(I-K)^{-1}R(t)^{-1}B(t)^{\top}P(t)]X(t)\\
&-[B(t)R(t)^{-1}B(t)^{\top}+(B(t)K+\widetilde{B}(t))(I-K)^{-1}\\
&\times R(t)^{-1}B(t)^{\top}]\psi(t)\}dt,\\
&d\psi(t)\\=&~-\{[A(t)^{\top}-P(t)B(t)R(t)^{-1}B(t)^{\top}\\
&-P(t)(B(t)K+\widetilde{B}(t))(I-K)^{-1}R(t)^{-1}B(t)^{\top}]\psi(t)\\
&-P(t)(B(t)K+\widetilde{B}(t))(I-K)^{-1}\\
&\times R(t)^{-1}B(t)^{\top}P(t)X(t)\}dt,\\
&X(0)=~x,~\psi(T)=0,
\end{aligned}
\right.
\end{equation}
and $P(\cdot)$ is the solution to \eqref{Theta}. Moreover, if there exist two nonnegative constants $\lambda_1$ and $\lambda_2$ with $\lambda_1+\lambda_2>0$, such that for any $x,y\in\mathbb{R}^n$, it holds that
\begin{equation}\label{MC}
\begin{aligned}
&\langle P(t)(B(t)K+\widetilde{B}(t))(I-K)^{-1}R(t)^{-1}B(t)^{\top}P(t)x,x\rangle\\
&-\langle B(t)R(t)^{-1}B(t)^{\top}y,y\rangle\\
&-\langle (B(t)K+\widetilde{B}(t))(I-K)^{-1} R(t)^{-1}B(t)^{\top}y,y  \rangle\\
\leq& -\lambda_1|x|^2-\lambda_2|y|^2,
\end{aligned}
\end{equation}
thus CC system \eqref{CC} admits a unique solution $(X(\cdot),\psi(\cdot))$.

\textbf{Proof.}
From \eqref{statefeedback}, it follows
\begin{equation}
\begin{aligned}\label{m1}
m(t)=&\underset{N\rightarrow\infty}{\lim}u^{*(N,-i)}(t)=\underset{N\rightarrow\infty}{\lim}\frac{1}{N-1}\underset{j=1,j\neq i}{\overset{N}\sum} u^{\ast}_j(t)\\
=&-R(t)^{-1}B(t)^{\top}\Big(P(t)\underset{N\rightarrow\infty}{\lim}\frac{1}{N-1}\underset{j=1,j\neq i}{\overset{N}\sum}\widehat{X}^{\ast}_{j}(t)\\
&+\psi(t)\Big)+Km(t),
\end{aligned}
\end{equation}
where $P(\cdot)$ and $\psi(\cdot)$ are the solutions to \eqref{Theta} and \eqref{psi}, respectively.

 Due to the setting of partial observation structure, the optimal state filtering $\widehat{X}^{*}_{i}(\cdot)$ is only driven by the new individual noise $\widehat{W}_{i}(\cdot)$, see \eqref{closedstate}. It shows from the analysis in Remark 2.3 that $\widehat{W}_{i}(\cdot)$ and thus $\widehat{X}^{\ast}_{i}(\cdot)$ are i.i.d., for $1\leq i\leq N.$ Consequently, by the law of large numbers, it implies
\begin{equation}\label{xlimit}
\underset{N\rightarrow\infty}{\lim} \frac{1}{N-1}\underset{j=1,j\neq i}{\overset{N}\sum}\widehat{X}^{\ast}_{j}(t)= \mathbb{E}[\widehat{X}^{\ast}_{i}(t)].
\end{equation}

Combining with \eqref{m1} and \eqref{xlimit}, we derive
\begin{equation}\label{m2}
m(t)=-(I-K)^{-1}R(t)^{-1}B(t)^{\top}(P(t)\mathbb{E}[\widehat{X}^{\ast}_{i}(t)]+\psi(t)).
\end{equation}

Substituting \eqref{m2} into \eqref{closedstate} and then taking expectation on both side of \eqref{closedstate}, we can get
\begin{equation}
\begin{aligned}\label{exi}
&d\mathbb{E}[\widehat{X}^{\ast}_{i}(t)]\\=&~\{[A(t)-B(t)R(t)^{-1}B(t)^{\top}P(t)\\
&-(B(t)K+\widetilde{B}(t))(I-K)^{-1}R(t)^{-1}B(t)^{\top}P(t)]\mathbb{E}[\widehat{X}^{\ast}_{i}(t)]\\
&-[B(t)R(t)^{-1}B(t)^{\top}+(B(t)K+\widetilde{B}(t))(I-K)^{-1}\\
&\times R(t)^{-1}B(t)^{\top}]\psi(t)\}dt.
\end{aligned}
\end{equation}

Therefore, $\mathbb{E}[\widehat{X}^{\ast}_{i}(\cdot)]=\mathbb{E}[\widehat{X}^{\ast}_{j}(\cdot)]$, for $i,j=1,2,\ldots,N$, which indicates that $\mathbb{E}[\widehat{X}^{\ast}_{i}(\cdot)]$ is independent of $i$. For simplicity, we let $X(\cdot)=\mathbb{E}[\widehat{X}^{\ast}_{i}(\cdot)]$. Recalling \eqref{m2}, we check that $m(\cdot)$ can be represented as \eqref{m}. Moreover, \eqref{exi} indeed deduces to the first equation in \eqref{CC}. Plugging \eqref{m} into \eqref{psi}, we obtain the updated equation for $\psi(\cdot)$ in \eqref{CC}.

The equation \eqref{CC} is called the NCE system or CC system of Problem (LPO), which is a coupled forward-backward ordinary differential equation. Using the results in \cite{Wu1997} or \cite{Peng1999}, we can verify that the monotonicity condition \eqref{MC} is sufficient for the well-posedness of \eqref{CC}.

By \eqref{m}, \eqref{CC} and the assumptions of the coefficients, we conclude that $m(\cdot)\in L^2(0,T;\mathbb{R}^k)$ is a deterministic function, which completes the proof.
 \hfill$\square$

\textbf{Remark 2.4.}
It is necessary to compare some works on partially observed MFGs \cite{Firoozi2020} and \cite{Bensoussan2019} with our paper. In \cite{Firoozi2020}, the authors studied a class of LQG major-minor MFGs with partial observations of all agents. Without considering major player, \cite{Bensoussan2019} studied a class of LQG MFGs with partial observation structure for individual agents. We note that the interactions between the individuals in both papers are achieved through the state average, and the classical separation principle is used. By contrast, we first study a control average LQ MFG with partial observation and then solve them by the backward separation principle in the current work.

\section{$\varepsilon$-Nash Equilibrium Analysis}
In this section, we will verify that the decentralized optimal control above-mentioned satisfies the $\varepsilon$-Nash equilibrium property. For this purpose, let us first recall the definition of $\varepsilon$-Nash equilibrium.

\textbf{Definition 3.1.} The strategy set $\boldsymbol{u}^{\ast}(\cdot)=(u^{\ast}_1(\cdot),\ldots,u^{\ast}_N(\cdot))$, where $u^{\ast}_i(\cdot)\in \mathcal{U}_i$, $1\leq i \leq N$, is called an $\varepsilon$-Nash equilibrium with respect to cost $\mathcal{J}_i$, if there exists an $\varepsilon=\varepsilon(N)$ and $\underset{N\rightarrow\infty}{\lim }\varepsilon(N)=0$, such that
\begin{equation*}
\mathcal{J}_{i}(u^{\ast}_{i}(\cdot ),\boldsymbol{u}^{\ast}_{-i}(\cdot ))\leq \mathcal{J}_{i}(u_{i}(\cdot ),\boldsymbol{u}^{\ast}_{-i}(\cdot ))+\varepsilon,
\end{equation*}
where $u_{i}(\cdot )\in \mathcal{U}_i$ is an alternative strategy applied by agent $i$.

With $u^{\ast}_i(\cdot)$ given by \eqref{statefeedback} and $\widehat{X}^{\ast}_{i}(\cdot)$ determined by \eqref{closedstate}, suppose $X^{\dag}_i(\cdot)$ is the centralized state satisfying the following equation
\begin{equation}
\left\{
\begin{aligned}
&dX^{\dag}_i(t)\\
=&~[A(t)X^{\dag}_i(t)-B(t)R(t)^{-1}B(t)^{\top}(P(t)\widehat{X}^{\ast}_{i}(t)\\
&+\psi(t))-\frac{1}{N-1}\underset{j=1,j\neq i}{\overset{N}\sum}\widetilde{B}(t)R(t)^{-1}B(t)^{\top}(P(t)\widehat{X}^{\ast}_{j}(t)\\
&+\psi(t))+(B(t)+\widetilde{B}(t))Km(t)]dt+\sigma
(t)dW_{i}(t),\\
&X^{\dag}_i(0)=x,
\end{aligned}
\right.
\end{equation}
and the observation process follows
\begin{equation}\label{centralizedobservation}
\left\{
\begin{aligned}
dV^{\dag}_{i}(t)=&~[F(t)X^{\dag}_{i}(t)+G(t)]dt+H(t)d\overline{W}_{i}(t), \\
V^{\dag}_{i}(0)=&~0.
\end{aligned}%
\right.
\end{equation}
Similarly, $X^{\ast}_i(\cdot)$ is the decentralized state satisfying
\begin{equation}
\left\{
\begin{aligned}\label{destate}
&dX^{\ast}_i(t)\\
=&~[A(t)X^{\ast}_i(t)-B(t)R(t)^{-1}B(t)^{\top}(P(t)\widehat{X}^{\ast}_{i}(t)\\
&+\psi(t))+(B(t)K+\widetilde{B}(t))m(t)]dt+\sigma(t)dW_{i}(t),\\
&X^{\ast}_i(0)=~x
\end{aligned}
\right.
\end{equation}
and the corresponding observation process $V_i^{\ast}(\cdot)$ still satisfies equation \eqref{centralizedobservation} with $X^{\dag}_i(\cdot)$ replaced by $X^{\ast}_i(\cdot)$.

Since $m(\cdot)$ has been presented in Proposition 2.2, the difference $X^{\dag}_i(\cdot)-X^{\ast}_i(\cdot)$ satisfies
\begin{equation}
\left\{
\begin{aligned}
&d(X^{\dag}_i(t)-X^{\ast}_i(t))=\Big\{A(t)(X^{\dag}_i(t)-X^{\ast}_i(t))-\widetilde{B}(t)R(t)^{-1}\\
&\qquad \qquad \times B(t)^{\top}P(t)\big[\frac{1}{N-1}\underset{j=1,j\neq i}{\overset{N}\sum}\widehat{X}^{\ast}_{j}(t)-X(t)\big]\Big\}dt,\\
&X^{\dag}_i(0)-X^{\ast}_i(0)=~0,
\end{aligned}
\right.
\end{equation}
where $X(\cdot)$ follows the forward  component of \eqref{CC}.

\textbf{Lemma 3.1.} Let (A1) and (A2) hold. It follows
\begin{equation}\label{1a}
\underset{1\leq i\leq N}{\sup }\Big[ \underset{0\leq t\leq T}{\sup }\mathbb{%
E}|X^{\dag}_{i}(t)-X^{\ast}_{i}(t)|^{2}\Big] =O\left( \frac{1}{N}\right),
\end{equation}
\begin{equation}\label{1b}
|\mathcal{J}_{i}\left( u^{\ast}_{i}(\cdot),\boldsymbol{u}^{\ast}_{-i}(\cdot)\right) -J_{i}\left(
u^{\ast}_{i}(\cdot)\right) |=O\left( \frac{1}{\sqrt{N}}\right).
\end{equation}

\textbf{Proof.} By the basic inequality and boundedness of coefficients, there exists a positive constant $C$ such that
\begin{equation}\label{xe}
\begin{aligned}
\mathbb{E}|X^{\dag}_i(t)-&X^{\ast}_i(t)|^{2}\leq C \int_0^t\mathbb{E}|X^{\dag}_i(s)-X^{\ast}_i(s)|^{2}ds\\
&+C\int_0^t\mathbb{E}\big|\frac{1}{N-1}\underset{j=1,j\neq i}{\overset{N}\sum}\widehat{X}^{\ast}_{j}(s)-X(s)\big|^2ds.
\end{aligned}
\end{equation}
From \eqref{closedstate} and \eqref{m}, we have
\begin{equation}
\left\{
\begin{aligned}\label{sumxhat}
&d\frac{1}{N-1}\underset{j=1,j\neq i}{\overset{N}\sum}\widehat{X}^{\ast}_{j}(t)\\
=&~\Big\{(A(t)-B(t)R(t)^{-1}B(t)^{\top}P(t))
\frac{1}{N-1}\underset{j=1,j\neq i}{\overset{N}\sum}\widehat{X}^{\ast}_{j}(t)\\
&~~-(B(t)K+\widetilde{B}(t))(I-K)^{-1}R(t)^{-1}B(t)^{\top}P(t)X(t)\\
&~~-[B(t)R(t)^{-1}B(t)^\top+(B(t)K+\widetilde{B}(t))(I-K)^{-1}\\
&~~\times R(t)^{-1}B(t)^{\top}]\psi(t)\Big\}dt\\
&~~+\frac{1}{N-1}\underset{j=1,j\neq i}{\overset{N}\sum}\Pi (t)F(t)^{\top }(H(t)^\top)^{-1}d\widehat{W}_{j}(t),\\
&\frac{1}{N-1}\underset{j=1,j\neq i}{\overset{N}\sum}\widehat{X}^{\ast}_{j}(0)=~x.
\end{aligned}
\right.
\end{equation}
\\
Combining \eqref{CC} with \eqref{sumxhat}, it follows
\begin{equation}
\left\{
\begin{aligned}
&d\Big(\frac{1}{N-1}\underset{j=1,j\neq i}{\overset{N}\sum}\widehat{X}^{\ast}_{j}(t)-X(t)\Big)
\\
&=~\Big[(A(t)-B(t)R(t)^{-1}B(t)^{\top} P(t))\\
&\qquad \qquad \quad \times \Big(\frac{1}{N-1}\underset{j=1,j\neq i}{\overset{N}\sum}\widehat{X}^{\ast}_{j}(t)-X(t)\Big)\Big]dt\\
&~~~~+\frac{1}{N-1}\underset{j=1,j\neq i}{\overset{N}\sum}\Pi (t)F(t)^{\top }(H(t)^\top)^{-1}d\widehat{W}_{j}(t),\\
&\frac{1}{N-1}\underset{j=1,j\neq i}{\overset{N}\sum}\widehat{X}^{\ast}_{j}(0)-X(0)=~0.
\end{aligned}
\right.
\end{equation}
\\
By integrating and taking expectation on both sides of the above equation, we have
\begin{equation*}
\begin{aligned}
&\mathbb{E}\big|\frac{1}{N-1}\underset{j=1,j\neq i}{\overset{N}\sum}\widehat{X}^{\ast}_{j}(t)-X(t)\big|^2\\
\leq & C\int_0^t\mathbb{E}\big|\frac{1}{N-1}\underset{j=1,j\neq i}{\overset{N}\sum}\widehat{X}^{\ast}_{j}(s)-X(s)\big|^2ds\\
&+C\int_0^t\mathbb{E}\big|\frac{1}{N-1}\underset{j=1,j\neq i}{\overset{N}\sum}\Pi (s)F(s)^{\top }(H(s)^\top)^{-1}\big|^2ds,\\
\leq & C\Big(1+\int_0^t\mathbb{E}\big|\frac{1}{N-1}\underset{j=1,j\neq i}{\overset{N}\sum}\widehat{X}^{\ast}_{j}(s)-X(s)\big|^2ds\Big),
\end{aligned}
\end{equation*}
which implies
\begin{equation}\label{sume}
\mathbb{E}\big|\frac{1}{N-1}\underset{j=1,j\neq i}{\overset{N}\sum}\widehat{X}^{\ast}_{j}(t)-X(t)\big|^2=O\left(\frac{1}{N}\right).
\end{equation}

From \eqref{xe}, \eqref{sume} and the Gronwall's inequality, we obtain the first assertion of this lemma.

By \eqref{cost} and \eqref{lcost}, the difference $\mathcal{J}_{i}\left( u^{\ast}_{i}(\cdot),\boldsymbol{u}^{\ast}_{-i}(\cdot)\right) -J_{i}\left(u^{\ast}_{i}(\cdot)\right)$ can be rewritten as
\begin{equation}
\begin{aligned}\label{dJ}
&\mathcal{J}_{i}\left( u^{\ast}_{i}(\cdot),\boldsymbol{u}^{\ast}_{-i}(\cdot)\right) -J_{i}\left(
u^{\ast}_{i}(\cdot)\right)\\
=&\frac{1}{2}\mathbb{E}\Big\{\int_0^T[\langle Q(t)X^{\dag}_i(t),X^{\dag}_i(t)\rangle-\langle Q(t)X^{\ast}_i(t),X^{\ast}_i(t)\rangle]dt\\
&+\langle MX^{\dag}_i(T),X^{\dag}_i(T)\rangle-\langle MX^{\ast}_i(T),X^{\ast}_i(T)\rangle\Big\}.
\end{aligned}
\end{equation}
\\
Noting the fact that $\underset{0\leq t\leq T}{\sup}\mathbb{E}|X^{\ast}_i(t)|^2<\infty$ and from \eqref{1a}, it holds
\begin{equation}
\begin{aligned}\label{firstJ}
&\Big|\mathbb{E}\int_0^T[\langle Q(t)X^{\dag}_i(t),X^{\dag}_i(t)\rangle-\langle Q(t)X^{\ast}_i(t),X^{\ast}_i(t)\rangle]dt\Big|\\
\leq&C\int_0^T\mathbb{E}|X^{\dag}_i(t)-X^{\ast}_i(t)|^2dt\\
&+C\int_0^T\mathbb{E}[|X^{\dag}_i(t)-X^{\ast}_i(t)||X^{\ast}_i(t)|]dt\\
\leq&C\int_0^T\mathbb{E}|X^{\dag}_i(t)-X^{\ast}_i(t)|^2dt\\
&+C\int_0^T(\mathbb{E}|X^{\dag}_i(t)-X^{\ast}_i(t)|^2)^{\frac{1}{2}}(\mathbb{E}|X^{\ast}_i(t)|^2)^{\frac{1}{2}}dt\\
=&O\left(\frac{1}{\sqrt{N}}\right),
\end{aligned}
\end{equation}
where the last inequality is due to the Cauchy-Schwartz inequality.

Similar techniques can be applied to the terminal term of \eqref{dJ}, which generates
\begin{equation}\label{secondJ}
|\mathbb{E}[\langle MX^{\dag}_i(T),X^{\dag}_i(T)\rangle-\langle MX^{\ast}_i(T),X^{\ast}_i(T)\rangle]|=O\left(\frac{1}{\sqrt{N}}\right).
\end{equation}

Based on \eqref{dJ}-\eqref{secondJ}, the desired results are obvious.  \hfill$\square$

In the following, we introduce some perturbations to Problem (LPO). Specifically, for fixed $i$, suppose the agent $i$ can take any $u_i(\cdot)\in\mathcal{U}_i$ and other agents $j$, $1\leq j \leq N, j\neq i$, still hold the decentralized optimal control $u^{\ast}_j(\cdot)$. Accordingly, the centralized trajectories of the agent $i$ and $j$ with perturbation read
\begin{equation}
\left\{
\begin{aligned}
&d\alpha_{i}(t)\\
=&~\Big[A(t)\alpha_{i}(t)+B(t)u_{i}(t)+\frac{1}{N-1}\underset{\kappa=1,\kappa\neq i}{\overset{N}\sum} \widetilde{B}(t)u_\kappa(t)\Big]dt \\
&+\sigma (t)dW_{i}(t), \\
&d\alpha_{j}(t)\\=&~\Big[A(t)\alpha_{j}(t)+B(t)u^{\ast}_{j}(t)+\frac{1}{N-1}\widetilde{B}(t)\Big(\underset{\kappa=1,\kappa\neq i,j}{\overset{N}\sum} u^{\ast}_\kappa(t)\\
&+u_i(t)\Big)\Big]dt +\sigma (t)dW_{j}(t), \\
&\alpha_{i}(0)=~x,~\alpha_{j}(0)=~x,
\end{aligned}
\right.
\end{equation}
and the perturbed observation processes satisfy
\begin{equation}
\left\{
\begin{aligned}\label{perturbedobservation}
dV^{\alpha}_{i}(t)=&~[F(t)\alpha_{i}(t)+G(t)]dt+H(t)d\overline{W}_{i}(t), \\
dV^{\alpha}_{j}(t)=&~[F(t)\alpha_{j}(t)+G(t)]dt+H(t)d\overline{W}_{j}(t), \\
V^{\alpha}_{i}(0)=&~0,~V^{\alpha}_{j}(0)=~0.
\end{aligned}%
\right.
\end{equation}
\\
Moreover, the auxiliary state of the agent $i$ satisfies
\begin{equation}
\left\{
\begin{aligned}
d\beta_{i}(t)=&~\Big[A(t)\beta_{i}(t)+B(t)u_{i}(t)+ \widetilde{B}(t)m(t)\Big]dt \\
&+\sigma (t)dW_{i}(t), \\
\beta_{i}(0)=&~x,
\end{aligned}
\right.
\end{equation}
and the corresponding observation process $V^{\beta}_{i}(\cdot)$ follows the first equation of \eqref{perturbedobservation} with $\alpha_{i}(\cdot)$ replaced by $\beta_{i}(\cdot)$.

Using a similar procedure as the previous one, we can derive the following result. For saving space, we do not prove it here.

\textbf{Lemma 3.2.} Let (A1) and (A2) hold. Then it follows
\begin{equation}
\underset{1\leq i\leq N}{\sup }\Big[ \underset{0\leq t\leq T}{\sup }\mathbb{%
E}|\alpha_{i}(t)-\beta_{i}(t)|^{2}\Big] =O\left( \frac{1}{N}\right),
\end{equation}
\begin{equation}
|\mathcal{J}_{i}\left( u_{i}(\cdot),\boldsymbol{u}^{\ast}_{-i}(\cdot)\right) -J_{i}\left(
u_{i}(\cdot)\right) |=O\left( \frac{1}{\sqrt{N}}\right).
\end{equation}
\\
Based on the results of Lemmas 3.1 and 3.2, we can get the main theorem of this section.

\textbf{Theorem 3.1.} Let (A1) and (A2) hold. Then  $\boldsymbol{u}^{\ast}(\cdot)=(u^{\ast}_1(\cdot),\ldots,u^{\ast}_N(\cdot))$, where $u^{\ast}_i(\cdot)$ is represented by \eqref{statefeedback}, is the $\varepsilon$-Nash equilibrium of Problem (LPO).

\textbf{Proof.} Using the results mentioned in Lemmas 3.1 and 3.2, we obtain
\begin{equation*}
\begin{aligned}
&\mathcal{J}_{i}\left( u^{\ast}_{i}(\cdot),\boldsymbol{u}^{\ast}_{-i}(\cdot)\right) =J_{i}\left(
u^{\ast}_{i}(\cdot)\right) +O\left( \frac{1}{\sqrt{N}}\right)\\
\leq &J_{i}\left( u_{i}(\cdot)\right) +O\left( \frac{1}{\sqrt{N}}\right)
=\mathcal{J}_{i}\left( u_{i}(\cdot),\boldsymbol{u}^{\ast}_{-i}(\cdot)\right) +O\left( \frac{1}{%
\sqrt{N}}\right),
\end{aligned}
\end{equation*}
thus the $\varepsilon$-Nash equilibrium property can be verified with $\varepsilon=O\left( \frac{1}{%
\sqrt{N}}\right)$.      \hfill$\square$

\section{An Illustrative Example}
In this section, we use the foregoing results to solve a cash management problem, which has important values in both theoretical and practical aspects. For simplicity, we consider a one-dimensional case, and the control average term only appears in the state equation.

For $1\leq i \leq N$, suppose that the liability process of the firm $i$ is controlled by
\begin{equation*}
-dL_i(t)=\Big[b(t)u_i(t)+\tilde{b}(t)u^{(N,-i)}(t)\Big]dt+\sigma(t)dW_i(t),
\end{equation*}
where $u_i(\cdot)$ denotes the control strategy of the firm $i$ and usually stands for the rate of capital injection or withdrawal. Moreover, the term $u^{(N,-i)}(\cdot):=\frac{1}{N-1}\sum_{j=1,j\neq i}^N u_j(\cdot)$ is introduced to minor the impact of others' decisions on the firm $i$. The diffusion term $\sigma(t)dW_i(t)$ denotes the individual liability risk, which may vary from firm to firm.

Assume the initial endowment of the firm $i$ is $x$, and it can only invest in a money account with compounded interest rate $a(\cdot)$. Thus, the cash-balance process of the firm $i$ is subject to
\begin{equation*}
\begin{aligned}
&X_i(t)\\
=&\exp\big\{\int_0^ta(s)ds\big\}\Big(x-\int_0^t\exp\big\{-\int_0^sa(\tau)d\tau\big\}dL_i(s)\Big),
\end{aligned}
\end{equation*}
which can be rewritten in the following differential form:
\begin{equation}\label{appstate}
\left\{
\begin{aligned}
dX_{i}(t)=&~[a(t)X_{i}(t)+b(t)u_{i}(t)+\tilde{b}(t)u^{(N,-i)}(t)]dt \\
&+\sigma (t)dW_{i}(t), \\
X_{i}(0)=&~x.
\end{aligned}
\right.
\end{equation}
\\
Due to technical limitations and the dispersion of account information, each company can only observes the cash-balance process through the corresponding stock price, which satisfies
\begin{equation*}
\left\{
\begin{aligned}
dS_i(t)=&~S_i(t)[f(t)X_i(t)+g(t)+\frac{1}{2}h(t)^2]dt\\
&+h(t)d\overline{W}_i(t),\\
S_i(0)=&~1.
\end{aligned}
\right.
\end{equation*}
\\
For $1\leq i \leq N$, let $V_i(\cdot)=\log S_i(\cdot)$, then
\begin{equation}
\left\{
\begin{aligned}\label{appobservation}
dV_i(t)=&~[f(t)X_i(t)+g(t)]dt+h(t)d\overline{W}_i(t),\\
V_i(0)=&~0.
\end{aligned}
\right.
\end{equation}
\\
Suppose that each firm wants to minimize the cost functional in form of
\begin{equation}\label{appcost}
\begin{aligned}
&\mathcal{J}_{i}(u_{i}(\cdot ),\boldsymbol{u}_{-i}(\cdot ))\\
=&\frac{1}{2}\mathbb{E}\Big\{\int_0^T(u_i(t)-r(t))^2dt+(X_i(T)-l)^2\Big\},
\end{aligned}
\end{equation}
where $r(\cdot)$ is an $\mathbb{R}$-valued deterministic and bounded function, $l$ is a constant.

The above cost functional consists of two terms: the first term reveals the derivation between the control $u_i(\cdot)$ and the benchmark $r(\cdot)$; and the second one measures the risk of terminal wealth. Then, the cash management problem with partial observation (CMP) can be formulated as follows.

\textbf{Problem(CMP)} Look for a capital injection or withdraw strategy set $\boldsymbol{u}^{\ast}(\cdot)=(u^{\ast}_1(\cdot),\ldots,u^{\ast}_N(\cdot))$ such that
\begin{equation*}
\mathcal{J}_{i}(u^{\ast}_{i}(\cdot ),\boldsymbol{u}^{\ast}_{-i}(\cdot ))=\underset{u_i(\cdot)\in\mathcal{U}_i}{\inf}
\mathcal{J}_{i}(u_{i}(\cdot ),\boldsymbol{u}^{\ast}_{-i}(\cdot )),~~ 1\leq i \leq N.
\end{equation*}
\\
Compared to \eqref{cost}, there exist linear terms  $r(t)u_i(t)$ and $lX_i(T)$ in \eqref{appcost}. But it brings no essential difficulties in our problem, which can be addressed similarly. Let us introduce $m(\cdot)$ and denote the limiting problem as Problem (LCMP), we can link it to the following stochastic Hamiltonian system
\begin{equation}
\left\{
\begin{aligned}\label{appH}
0=&~u^{\ast}_i(t)-r(t)+b(t)\mathbb{E}[Y^{\ast}_{i}(t)|\mathcal{F}^{%
V^{\ast}_{i}}_t],\\
dX^{\ast}_{i}(t)=&~[a(t)X^{\ast}_{i}(t)+b(t)u^{\ast}_{i}(t)+\tilde{b}(t)m(t)]dt\\
&+\sigma
(t)dW_{i}(t), \\
dY^{\ast}_{i}(t)=&-a(t)Y^{\ast}_{i}(t)dt+Z^{\ast}%
_{i}(t)dW_{i}(t), \\
X^{\ast}_{i}(0)=&~x,Y^{\ast}_{i}(T)=X^{\ast}_{i}(T)-l,%
\end{aligned}
\right.
\end{equation}
where $V_i^*(\cdot)$ satisfies \eqref{appobservation} in form with $X_i^*(\cdot)$ given by the second equation of \eqref{appH}.

Recall the notation $\widehat{\theta}_{i}(\cdot )=\mathbb{E}[\theta_{i}(\cdot )|\mathcal{F}^{%
V^{\ast}_{i}}_\cdot]$, $1\leq i\leq N$, then the optimal filtering equation reads
\begin{equation}
\left\{
\begin{aligned}\label{appfilter}
d\widehat{X}^{\ast}_{i}(t)=&~[a(t)\widehat{X}^{\ast}_{i}(t)-b(t)^2\widehat{Y}^{\ast}_{i}(t)+b(t)r(t)\\
&+\tilde{b}(t)m(t)]dt+\frac{\Pi (t)f(t)}{h(t)}d\widehat{W}_{i}(t),
\\
d\widehat{Y}^{\ast}_{i}(t)=&-a(t)\widehat{Y}^{\ast}_{i}(t)dt\\
&+\frac{f(t)(\widehat{Y^{\ast}_i(t)X^{\ast}_i(t)}-\widehat{Y}^{\ast}_i(t)\widehat{X}^{\ast}_i(t))}{h(t)}d\widehat{W}_{i}(t), \\
\widehat{X}^{\ast}_{i}(0)=&~x,~\widehat{Y}^{\ast}_{i}(T)=~\widehat{X}^{\ast}_{i}(T)-l,
\end{aligned}%
\right.
\end{equation}
with
\begin{equation*}
\widehat{W}_i(t)=\int_0^t \frac{f(s)(X^{\ast}_i(s)-\widehat{X}^{\ast}
_{i}(s))}{h(s)}ds+\overline{W}_i(t)
\end{equation*}
and $\Pi (t)=\mathbb{E}[(X^{\ast}_{i}(t)-\widehat{X}^{\ast}_{i}(t))^2|\mathcal{F}^{
V^{\ast}_{i}}_t]$ is the solution to
\begin{equation*}
\left\{
\begin{aligned}
&\dot{\Pi}(t)-2a(t)\Pi (t)+\Big(\frac{\Pi
(t)f(t)}{h(t)}\Big)^2-\sigma (t)^2=0, \\
&\Pi (0)=0.%
\end{aligned}%
\right.
\end{equation*}
\\
From \eqref{appH}, we can derive the following open-loop decentralized optimal control
\begin{equation}\label{appopen}
u^{\ast}_i(t)=-b(t)\widehat{Y}^{\ast}_{i}(t)+r(t),
\end{equation}
where $\widehat{Y}^{\ast}_{i}(\cdot)$ follows the second equation of \eqref{appfilter}.

Furthermore, we introduce the following equations:
\begin{equation*}
\left\{
\begin{aligned}
&\dot{P}(t)+2a(t)P(t)-b(t)^2 P(t)^2=0,\\
&P(T)=1,
\end{aligned}%
\right.
\end{equation*}
and
\begin{equation*}
\left\{
\begin{aligned}
&\dot{\psi}(t)+(a(t)-b(t)^2P(t))\psi(t)\\
&\quad  +P(t)b(t)r(t)+P(t)\tilde{b}(t)m(t)=0,\\
&\psi(T)=-l,
\end{aligned}%
\right.
\end{equation*}
then the feedback representation of decentralized optimal control satisfies
\begin{equation}\label{appclosed}
u^{\ast}_i(t)=-b(t)(P(t)\widehat{X}^{\ast}_{i}(t)+\psi(t))+r(t),
\end{equation}
with
\begin{equation*}
\left\{
\begin{aligned}
d\widehat{X}^{\ast}_{i}(t)=&~[(a(t)-b(t)^2P(t))\widehat{X}^{\ast}_{i}(t)-b(t)^2\psi(t)\\
&+b(t)r(t)+\tilde{b}(t)m(t)]dt+\frac{\Pi (t)f(t)}{h(t)}d\widehat{W}_{i}(t),
\\
\widehat{X}^{\ast}_{i}(0)=&~x.
\end{aligned}%
\right.
\end{equation*}

By the results in Proposition 2.2, we can derive the similar CC system with solution $(X(\cdot),\psi(\cdot))$.  We further let $\psi(t)=\Gamma(t)X(t)+\Lambda(t)$ to decouple the related coupled CC system, where $\Gamma(\cdot)$ and $\Lambda(\cdot)$ are governed by
\begin{equation*}
\left\{
\begin{aligned}
&\dot{\Gamma}(t)+2[a(t)-P(t)(b(t)+\tilde{b}(t))b(t)]\Gamma(t)\\
&\qquad-(b(t)+\tilde{b}(t))b(t)\Gamma(t)^2-b(t)\tilde{b}(t)P(t)^2=0,\\
&\Gamma(T)=0,
\end{aligned}%
\right.
\end{equation*}
\begin{equation*}
\left\{
\begin{aligned}
&\dot{\Lambda}(t)+[a(t)-(P(t)+\Gamma(t))(b(t)+\tilde{b}(t))b(t)]\Lambda(t)\\
&+(P(t)+\Gamma(t)(b(t)+\tilde{b}(t))r(t)=0,\\
&\Lambda(T)=-l.
\end{aligned}%
\right.
\end{equation*}

In such a way, we can first solve $P(\cdot)$ then $\Gamma(\cdot)$. Once $P(\cdot)$ and $\Gamma(\cdot)$ are available, $\Lambda(\cdot)$ becomes an ordinary differential equation. After simple calculations, then the control average limit can be solved explicitly as follows:
\begin{equation*}
\begin{aligned}
&m(t)\\=&-b(t)\Big[(P(t)+\Gamma(t))\exp\big\{\int_0^t[a(s)-(P(s)+\Gamma(s))\\
&\times(b(s)+\tilde{b}(s))b(s)]ds \big\}\Big\{x-\int_0^t
(b(s)+\tilde{b}(s))(b(s)\Lambda(s)\\
&-r(s))\exp\big\{-\int_0^s[a(\tau)-(P(\tau)+\Gamma(\tau))(b(\tau)\\
&+\tilde{b}(\tau))b(\tau)]d\tau \big\}ds\Big\}+\Lambda(t)\Big]+r(t).
\end{aligned}
\end{equation*}
Based on the above analysis and theoretical results in Section 3, we obtain the following results.

\textbf{Proposition 4.1.} Problem (CPM) admits the $\varepsilon$-Nash equilibrium  $\boldsymbol{u}^{\ast}(\cdot)=(u^{\ast}_1(\cdot),\ldots,u^{\ast}_N(\cdot))$ with $u^{\ast}_i(\cdot)$ given by \eqref{appclosed}.

In order to visualize our theoretical results, some numerical simulations are provided. We consider a large-population system with $100$ individual firms. For simplicity, we let $a=0.5$, $b=0.2$, $\tilde{b}=0.5$, $ \sigma=10$, $f=2.8$, $g=6$, $h=4$, $r=15$, $l=3$, $T=10$, $x=3.5$.

As shown in the following, Fig. 1 and Fig. 2 present the numerical solutions of the corresponding Riccati equations. The sample paths for the decentralized optimal states filtering and decentralized optimal controls of $100$ individual firms are given in  Fig. 3 and Fig. 4. Moreover, Fig. 5 describes the error estimate between the control average $u^{*(100,-1)}$ and the corresponding limit $m$ with the unit values scales. The gap between the centralized and the decentralized
states in square integrable sense is demonstrated in Fig. 6, where the agent number $N$ grows from 1 to 100. The attached simulations illustrate the good performance of our theoretical results.

\begin{figure}[H]
\label{P}
\begin{center}
\includegraphics[height=6.27cm, width=3.48in]{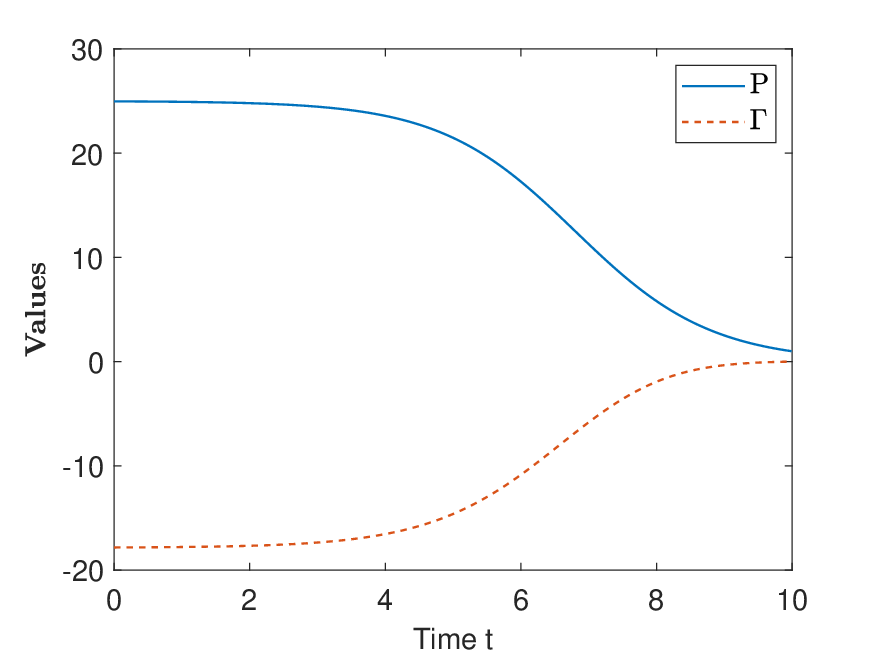}
% The printed column
% width is 8.4 cm.
% Size the figures
\end{center}
\par
\caption{The numerical solutions of $P$ and $\Gamma$}
% accordingly.
\end{figure}

\begin{figure}[H]
\label{Pii}
\begin{center}
\includegraphics[height=6.27cm, width=3.48in]{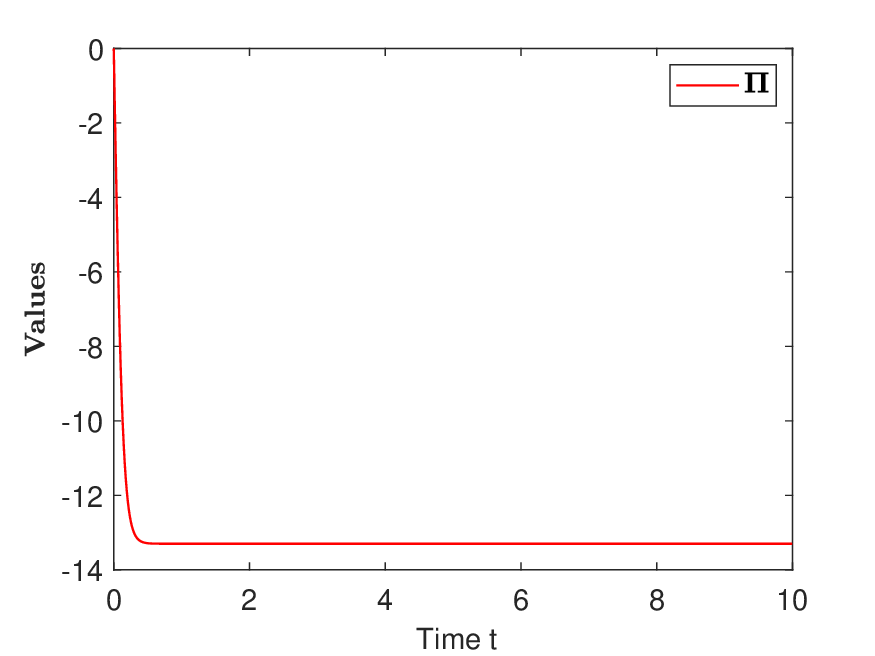}
% The printed column
% width is 8.4 cm.
% Size the figures
\end{center}
\par
\caption{The numerical solution of $\Pi$}
% accordingly.
\end{figure}

\begin{figure}[H]
\label{figstate}
\begin{center}
\includegraphics[height=6.27cm, width=3.48in]{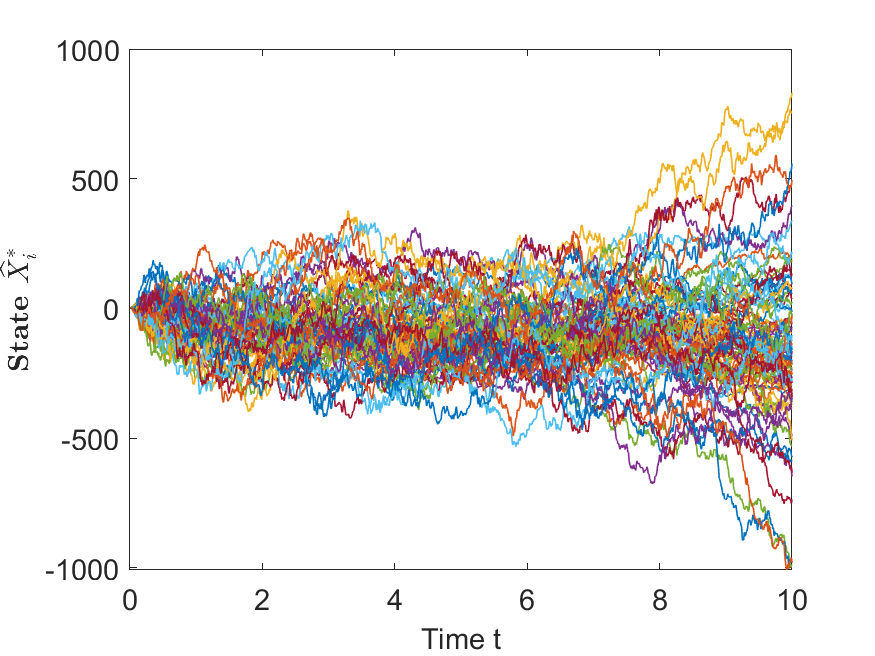}
% The printed column
% width is 8.4 cm.
% Size the figures
\end{center}
\par
\caption{The decentralized optimal states filtering}
% accordingly.
\end{figure}

\begin{figure}[H]
\label{figcontrol}
\begin{center}
\includegraphics[height=6.27cm, width=3.48in]{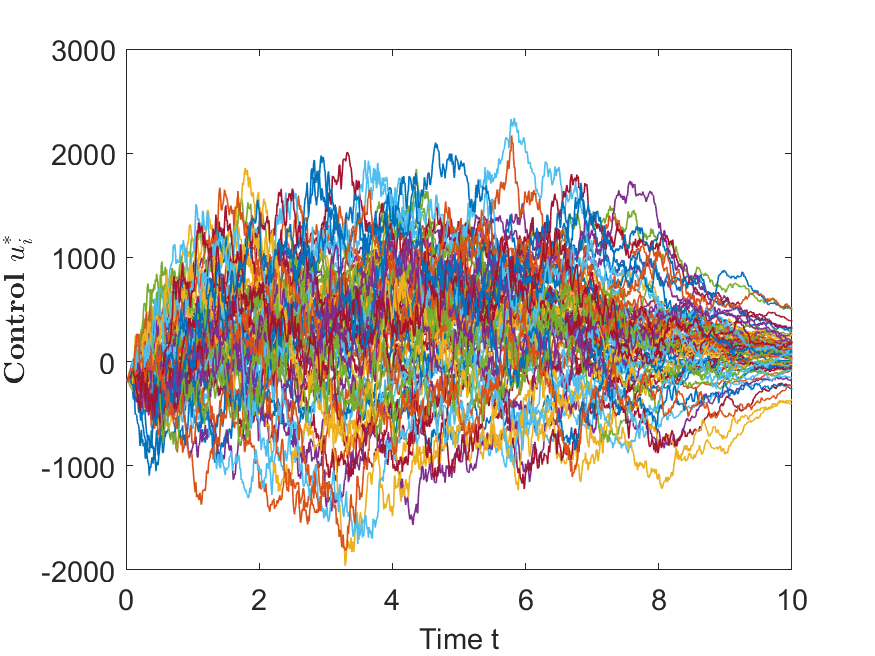}
% The printed column
% width is 8.4 cm.
% Size the figures
\end{center}
\par
\caption{The decentralized optimal controls}
% accordingly.
\end{figure}

\begin{figure}[H]
\label{figerror}
\begin{center}
\includegraphics[height=6.27cm, width=3.48in]{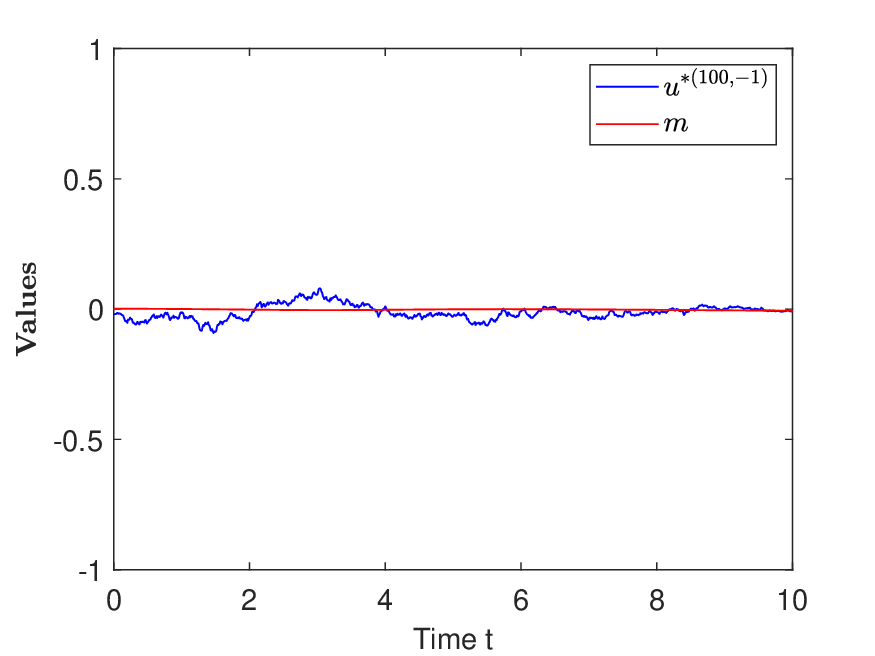}
% The printed column
% width is 8.4 cm.
% Size the figures
\end{center}
\par
\caption{The approximate error estimate }
% accordingly.
\end{figure}

\begin{figure}[H]
\label{figgap}
\begin{center}
\includegraphics[height=6.27cm, width=3.48in]{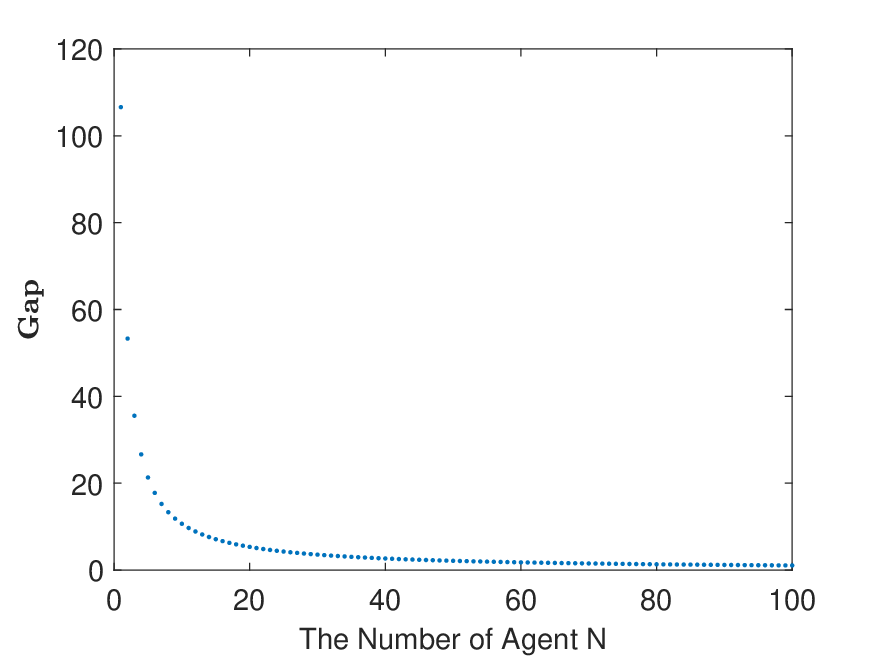}
% The printed column
% width is 8.4 cm.
% Size the figures
\end{center}
\par
\caption{The gap between centralized and decentralized
states}
% accordingly.
\end{figure}

\section{Conclusions}
This paper studies a class of LQ stochastic large-population problems with partial observation. The individual state process involved in our setting is not observed, where only a noisy observation related to it is gained. Using the MFG approach and the backward separation principle with a state decomposition technique, the decentralized optimal control can be obtained in the open-loop form through an FBSDE with the conditional expectation. The optimal filtering equation is also provided. By the decoupling method, decentralized optimal controls are presented as the feedback of states filtering. The CC system is discussed. Moreover, the related $\varepsilon$-Nash equilibrium property is verified. A cash management problem with partial observation is studied as an application in the end.

\end{document}